\documentclass[11pt,a4paper]{amsart}
\usepackage{amssymb}
\usepackage[backref=true]{biblatex}
\usepackage{color}
\usepackage[usenames,dvipsnames,svgnames,table]{xcolor}
\usepackage[colorlinks=true,
            linkcolor=red,
            citecolor=blue]{hyperref}

\usepackage{graphicx,amssymb,amsfonts,epsfig,amsthm,a4,amsmath,url}
\usepackage[latin1]{inputenc}
\usepackage{calligra}
\usepackage[T1]{fontenc}
\usepackage{amsmath, amsthm, amssymb, verbatim}
\usepackage{graphicx,amssymb,amsfonts,epsfig,amsthm,a4,amsmath,url,graphicx,enumerate}
\usepackage[latin1]{inputenc}
\usepackage{mathtools}
\usepackage{textcomp}
\usepackage{textpos}
\usepackage{tikz}

\newtheorem{thm}{Theorem}[section]
\newtheorem{cor}[thm]{Corollary}
\newtheorem{lem}[thm]{Lemma}

\newtheorem{prop}[thm]{Proposition}
\theoremstyle{definition}
\newtheorem{defn}[thm]{Definition}

\theoremstyle{remark}
\newtheorem{rem}[thm]{Remark}

\newtheorem{ex}[thm]{Example}

\numberwithin{equation}{section}

\newcommand{\II}{\mathcal{I}}

\newcommand{\LL}{\mathbf{L}}

\newcommand\restr[2]{{
  \left.\kern-\nulldelimiterspace 
  #1 
  \littletaller 
  \right|_{#2} 
  }}

\newcommand{\littletaller}{\mathchoice{\vphantom{\big|}}{}{}{}}

\begin{document}

\title[IME and Gromov's conjecture]{Integrable Measure Equivalence and Gromov's conjecture on nilpotent groups}

\author[Das]{ Kajal Das  }
\address{Woxsen University, Hyderabad, India}
\email{kdas.math@gmail.com}

\maketitle

\textbf{Abstract:}  In this article we prove a well-known conjecture  due to M. Gromov that  if two infinite finitely generated  torsion free nilpotent groups $\Gamma$ and $\Lambda$ are Quasi-isometric, then their Mal'cev Completions  are isomorphic as Lie groups. We prove this result in more general set-up. We prove that if two infinite finitely generated torsion free nilpotent groups  $\Gamma$ and $\Lambda$ are Integrable Measure Equivalent, then the Mal'cev completions of $\Gamma$ and $\Lambda$ are isomorphic as Lie groups. Equivalently, using Y. Shalom's result we prove that two infinite countable torsion free nilpotent groups  $\Gamma$ and $\Lambda$ are Integrable Measure Equivalent if and only if $\Gamma$ and $\Lambda$ are quasi-isometric if and only if the Mal'cev Completions of $\Gamma$ and $\Lambda$ are isomorphic as Lie groups. Moreover, we discuss some applications of our main theorem to box spaces of nilpotent groups, lifting of quasi-isometric embedding from $\Gamma$ to $\Lambda$ to a homo-morphism between their Mal'cev completions, and betti numbers of nilpotent groups. 
\vspace{5mm}


\vspace{5mm}

\textbf{Key terms:} Nilpotent groups, Nilpotent Lie groups, Mal'cev completion, Carnot groups, Carnot-Cratheodory Metric, Quasi-isometry, Topological Coupling, Integrable Measure Equivalence

\tableofcontents

\section{Introduction}

The class countable nilpotent groups or Nilpotent Lie groups is an interesting and mysterious class of groups since the inception of infinite group theory as this class is not classifiable. However, the \textit{large scale geometry} of nilpotent groups was alternative for studying such groups. The large scale geometry of nilpotent groups were started by M. Gromov by proving  a remarkable result in \cite{Gro81} that a finitely generated group of \textit{polynomial growth} is \textit{nilpotent} ( the result that a nilpotent has polynomial growth was already proved by H. Bass in \cite{Bas72} ). We refer to Subsection \ref{caygraph} for more details. Later,  the the large scale geometry of such groups were extensively studied  by P. Pansu in \cite{Pan83}, \cite{Pan89}. 

\vspace{2mm}

\textbf{Pansu's Theorem:}(\cite{Pan89} (Th\'{e}or\`{e}me 3, p.6)\textit{ Let $N$ be a finitely generated virtually nilpotent group equipped with some finite generating set $T$ and a left-invariant word-metric $d$. Then $(N, d /n, e)$ converges in pointed Gromov-Hausdorff topology to  to its \textit{asympototic cone}, which is a simply connected nilpotent Lie group $N_\mathbb{R}$ equipped with a left-invariant Carnot-Caratheodory (sub-Finsler) metric $d_{cc}$. $N_\mathbb{R}$ is called Carnot group associated with $N$. Moreover,  if two such groups $\Gamma$ and $\Lambda$ ( with symmetric generating sets $S_\Gamma$ and $S_\Lambda$, respectively) have bi-Lipschitz asymptotic cones (e.g., when are $\Gamma$ and $\Lambda$ are \textit{quasi-isometric}), then $\Gamma_\mathbb{R}$ and $\Lambda_{\mathbb{R}}$ isomorphic. }

\vspace{2mm} 

However, if  $\Gamma$ is a  finitely generated torsion free nilpotent groups,  A.I. Ma'cev proves in \cite{Mal51} that $\Gamma$ is isomorphic to a co-compact lattice in a simply connected nilpotent Lie group , denoted by $N(\Gamma)$, which is called the \textit{Mal'cev completion} of $\Gamma$.  We refer  to Subsection 2.5 for the definitions of  Mal'cev completion and Carnot group. 

The following is a well-known conjecture due to M. Gromov in geometric group theory which was known for long-time. 
 
  \vspace{2mm} 

\textbf{Gromov's Conjecture:} \textit{\textbf{If two infinite countable torsion free nilpotent groups $\Gamma$ and $\Lambda$ are quasi-isometric, then their Mal'cev Completions  are isomorphic as groups.}}

\vspace{2mm} 
 
 However, in his seminal article \cite{Gro93}, M. Gromov introduces \textit{Measure Equivalence} (ME) between two groups  as a measure theoretic counterpart of Quasi-isometry using  measurable dynamics of groups. In general, measure equivalence does not preserve large scale geometry of a group unlike quasi-isometry. But, if the`cocycles' arising from a measure equivalence relation satisfy an \textit{$\LL^p$-integrability condition}, it becomes an effective tool for studying lattices in Lie groups as observed by G. Margulis in \cite{Mar79} and Y. Shalom in \cite{Sha00}, \cite{Sha00}). The relation of Measure Equivalence together with $\LL^p$-integrability condition on the associated cocycles is called \textit{$\LL^p$-Measure Equivalence ($\LL^p$-ME)}. This new relation indeed preserves some large scale geometric properties of a finitely genearated group. The weakest form is when $p=1$, in which case the groups are called \textit{Integrable Measure Equivalent} (IME), and the strongest form is  when $p=\infty$, in which case  the groups are called \textit{Uniformly Measure Equivalent} (UME).There is an emerging interest in studying $\LL^p$-ME  in recent years both for amenable groups and non-amenable groups (see \cite{Sha04},\cite{BFS13},\cite{Aus16},\cite{DT16},\cite{Aus16'},\cite{Bow17},\cite{Das18}). 
 
\vspace{2mm} 
  
 In this article, our goal is to prove Gromov's conjecture in more general setting of IME . We prove the following  theorem: 

\begin{thm}\label{MainTheorem}(\textit{Main Theorem})
If two infinite countable torsion free nilpotent groups  $\Gamma$ and $\Lambda$ are Integrable Measure Equivalent, then the Mal'cev completions of $\Gamma$ and $\Lambda$ are isomorphic as Lie groups. 
\end{thm}

It has been proved by Y. Shalom in  \cite{Sha04} that \textit{if two amenable groups are quasi-isometric, then they are UME} . Since nilpotent groups are amenable and UME implies IME for two groups, our theorem proves Gromov's conjecture. Moreover, since any two co-compact lattices in a locally compact group are UME, we obtain that 

\begin{thm}\label{cormainthm}
The followings are equivalent: 

\begin{itemize}

\item  Two infinite countable torsion free nilpotent groups  $\Gamma$ and $\Lambda$ are Integrable Measure Equivalent; 

\item Two infinite countable torsion free nilpotent groups  $\Gamma$ and $\Lambda$ are Quasi-isometric; 

\item  The Mal'cev completions of $\Gamma$ and $\Lambda$ are isomorphic as Lie groups. 

\end{itemize}

\end{thm}

\section{Preliminaries} 

Since we write this article both for readers in ergodic theory and geometric group theory, we briefly recall  some definitions in both areas. 

\subsection{Cayley graph and growth of balls}\label{caygraph}

\begin{defn}(\textit{Cayley graph})
If $G$ is a finitely generated group generated by a finite symmetric set $S$, then we define the \textit{Cayley graph} $Cay(G,S)$ as follows: the set $V(G)$ of vertices are the group elements and two groups elements $g$ and $g'$ are connected by an edge if $g'=gs$, where $s\in S$. 
\end{defn}

We denote the $R$-radius ball of $G$ around the identity element $1$ by $B^G(1;R)$ for all $R\geq 0$. The cardinality $|B^G(1;R)|$ will be called the \textit{volume} of  $B^G(1;R)$ and $\alpha: [0,\infty)\rightarrow [0,\infty)$ will be called 
the \textit{growth function} of the Cayley graph $Cay(G,S)$. We say that a growth function $\alpha_2$ \textit{dominates} $\alpha_1$, denoted by $ \alpha_1 \prec \alpha_2$, if there exists a constant $\lambda\geq 1$ such that $\alpha_1(R)\leq \alpha(\lambda R)$ for all $R>0$. Two growth functions $\alpha_1$ and $\alpha_2$ are \textit{equivalent}, denoted by $\alpha_1 \sim \alpha_2$,  if $ \alpha_1 \prec \alpha_2$ and $ \alpha_2 \prec \alpha_1$. We note that $G$ is  generated by another 
symmetric generating set $T$ and $\alpha_1$ (resp. $\alpha_2$) is the growth function of the Cayley graph $Cay(G,S)$ (resp. $Cay(G,T)$), then $\alpha_1 \sim \alpha_2$.

\begin{rem}\label{grwothfnnp}
By H. Bass , a Cayley graph of a nilpotent  finitely generated group $G$ has growth function, say $\alpha$,  exactly polynomial of degree a of a certain positive integer $d(G)$ which can be characterized algebraically , i.e., $\alpha (R) \sim R^{d(G)}$ (see \cite{Bas72}). Precisely, $d(G)$ is equal to the \textit{Rank} of $G$, where $\mbox{Rank}(G)$ is defined as follows: Let $G_0\subset G_1\dots\subset G_k=G$ be the central series of $G$. Then, 

 $$ \mbox{Rank}(G):=\sum_{i=1}^k \mbox{Rank}\big( G_i/G_{i-1}\big).$$ 
\end{rem}

\subsection{Coarse equivalence and Quasi-isometry}\label{ceqi}

\begin{defn}(\textit{Coarse equivalence and Quasi-isometry for metric spaces})

Two metric spaces $X$ and $Y$  are \textit{Coarsely equivalent} if there exists a function $f$ from $X$ to $Y$ such that  for some monotonically increasing functions $\rho_{+}, \rho_{-}:[0,\infty)\rightarrow[0,\infty)$ going to $\infty$ and $D\geq 0$ such that 
 
\begin{equation}\label{cond1}
 \rho_{-}\big(d_X(x_1,x_2)\big) \leq d_Y\big(f(x_1),f(x_2)\big)\leq  \rho_{+}\big(d_X(x_1,x_2)\big) 
\end{equation} 
for all $x_1, x_2\in X$ and the image of $f$ is $D$-dense in $Y$, where $d_X$ and $d_Y$ are the distances in  $X$ and $Y$ respectively. If the first condition [\ref{cond1}] is only satisfied, $f$ is called a \textit{coarse embedding}. If $\rho_{+}(r)=Kr+C$ and $\rho_{-}(r)=\frac{r}{K}-C$ for some $K\geq 1$ and $C\geq 0$, we say that $X$ and $Y$ are $(K,C,D)$-\textit{Quasi-isometric} (QI). 

\end{defn}

\vspace{2mm}

It is well known that Coarse Equivalence can be defined equivalently by the following conditions: 

\begin{defn}\label{cecondn}(\cite{Pan89}, p.37)
Two metric space $X$ and $Y$ are \textit{coarsely equivalent} if and only if there exists two mappings $\phi:X\rightarrow Y$ 
and $\psi: Y\rightarrow X$ and monotonically increasing functions $\rho_{1}, \rho_{2}:[0,\infty)\rightarrow [0,\infty)$, $D\geq 0$ such that 

$$ d_Y\big(\phi(x_1), \phi(x_2)\big)\leq \rho_{1}\big(d_X(x_1,x_2)\big),\hspace*{1cm}d_X\big(\psi(y_1), \psi(y_2)\big)\leq \rho_{2}\big(d_Y(y_1,y_2)\big)$$

$$ d_X\big((\psi\circ\phi)(x),x\big)\leq C, \hspace*{1cm} d_Y\big((\phi\circ\psi)(y),y\big)\leq C $$
 for all $x_1, x_2, x\in X$ and $y_1,y_2,y\in Y$. 
\end{defn}

\begin{defn}(\textit{Coarse equivalence, Quasi-isometry for groups})
 Two countable groups $\Gamma$ (with a finite generating set  $S_\Gamma$) and $\Lambda$ (with a finite generating  $S_\Lambda$) are said to be \textit{Coarsely equivalent} (resp.  \textit{Quasi-isometric} ) if the Cayley graphs $(\Gamma, S_\Gamma)$ and  $(\Lambda, S_\Lambda)$ are Coarsely equivalent ( resp. Quasi-isometric) as metric spaces. 
 \end{defn}
 
 \begin{ex}
 The examples of quasi-isometric groups are as follows: Any two co-compact lattices in a locally compact  group (in particular, in a Lie group) are quasi-isometric ( we recall that a co-compact lattice is  a discrete subgroup such that the quotient of the ambient group by the discrete group is compact). 
 \end{ex}
 
M. Gromov proves the following theorem in his seminal article \cite{Gro93}: 
 
\begin{thm}\label{topocoupling}
The following are equivalent: 
\begin{itemize}
\item[1.] Two finitely generated groups $\Gamma$ and $\Lambda$ are Quasi-isometric; 
\item[2.] $\Gamma$ and $\Lambda$ have continuous , commuting, co-compact properly discontinuous action on a locally compact  space $X$.
\item[3.] $\Gamma$ and $\Lambda$ are Coarsely Equivalent. 
\end{itemize}
\end{thm}

  The locally compact space $X$ in the above-mentioned  is called \textit{Topological coupling} between the groups $\Gamma$ and $\Lambda$.  For example, if we take two co-compact lattices $\Gamma$ and $\Lambda$ in a same locally compact second countable Hausdorff group and we take the left action of $\Gamma$ on $G$ and the right action of $\Lambda$ on $G$, then they satisfy the above-mentioned `topological dynamical' relation. 
  
  \begin{rem}\label{qigrowth}
 It is easy to see that if two finitely generated nilpotent groups $\Gamma$ and $\Lambda$ are QI and $\Gamma$ (resp. $\Lambda$) has growth function $\alpha(R)\sim R^{d(\Gamma)} $ (resp. $\beta(R)\sim R^{d(\Lambda)}$), then 
 $d(\Gamma)=d(\Lambda)$ (see \cite{Gro81},  p. 57). In \cite{Aus16}, it is proved that if two nilpotent groups are IME, then  $d(\Gamma)=d(\Lambda)$. 
  \end{rem}
  
 \subsection{Gromov-Haudorff convergence}\label{ghconv}
 
 Suppose $X$ and $Y$ are two metric spaces and $f: X\rightarrow Y$ is a map. We define the `distortion' of $f$ by the 
following quantity: $$dis f := sup_{\{x_1,x_2\}} | d_Y(f(x_1),f(x_2))-d_X(x_1,x_2) |.$$ $f:X\rightarrow Y$ is called an $\epsilon$-isometry for 
some $\epsilon\geq0$ if $dis f \leq \epsilon$ and $Y$ is an $\epsilon$-neighborhood of $f(X)$. We sometimes say that `$f(X)$ is $\epsilon$-dense in 
$Y$' if $Y$ is an $\epsilon$-neighborhood of $f(X)$. There are several equivalent formulations of Gromov-Hausdorff 
convergence of compact metric spaces, we choose the following one for our purpose.  

\begin{defn}(\textit{Gromov-Hausdorff convergence}, \cite{BBI01}, p. 260)
A sequence $\{X_n\}$ of compact metric spaces 
converges to a compact metric space $X$ if there is a sequence $\{\epsilon_n\}_{n\in\mathbb{N}}$ of positive numbers and a sequence of maps 
$f_n: X_n \rightarrow X$ such that every $f_n$ is an $\epsilon_n$-isometry and $\epsilon_n\rightarrow 0$ as $n\rightarrow\infty$.
\end{defn}

\begin{lem}\label{relcomp}(Gromov's compactness criterion, \cite{Gro81}, p. 64)
A sequence $\{X_n\}_{n\in\mathbb{N}}$ of compact metric spaces converges to a compact metric space  ( or it is `relatively compact' in Gromov-Hausdorff metric) if only only if it satisfies the following criteria:
 
\begin{itemize}
 \item  there is a constant $E$ such that $\mbox{diam}\hspace*{1mm} X_n \leq E$ for all $n$;
 \item  for all $\epsilon>0$ there exists a natural number $N=N(\epsilon)$ such that
every $X_n$ contains an $\epsilon$-net consisting of at most $N$ points.
\end{itemize}

\end{lem}

\subsection{Measure Equivalence (ME) and $\LL^p$-ME}\label{melpme}

 As a measure theoretic counterpart of Quasi-isometry, M. Gromov introduces Measure Equivalence (ME) between two groups in the same article \cite{Gro93} using  measurable dynamics of groups. 
 
 \begin{defn}(\textit{Measure Equivalence(ME)})
 Two countable discrete groups $\Gamma$ and $\Lambda$ are called \textit{ Measure Equivalent}(ME) if there is a  standard Borel space $(X,\mu)$ (with a non-zero $\sigma$-finite measure $\mu$)  which admits commuting $\mu$-preserving actions of  $\Gamma$ and $\Lambda$  with finite-measure fundamental domains $X_{\Gamma}$ and $X_{\Lambda}$ respectively, i.e., 
$$X=\sqcup_{\gamma\in\Gamma}\hspace*{1mm}\gamma X_\Gamma=\sqcup_{\lambda\in\Lambda} \hspace*{1mm}\lambda X_\Lambda.$$ 
$(X,\mu)$  is called \textit{Measurable coupling space} between $\Gamma$ and $\Lambda$. 
\end{defn}

\begin{ex}
For example, if we take two lattices $\Gamma$ and $\Lambda$  in a locally compact second countable Hausdorff group $G$ with Haar Measure $\mu$  and we take the left action of $\Gamma$ on $G$ and the right action of $\Lambda$ on $G$, then they satisfy the above-mentioned measurable dynamical relation. 
\end{ex}

A Measure Equivalence relation between groups $\Gamma$ and $\Lambda$ gives rise to \textit{cocycle maps}  $\alpha:\Gamma\times X_\Lambda\rightarrow \Lambda$  (resp.\ $\beta:\Lambda\times X_\Gamma\rightarrow \Gamma$)  defined by the following rule: for all $x\in X_{\Lambda}$, and all $\gamma\in \Gamma$, $\alpha(\gamma, x)\gamma x \in X_\Lambda$ (and similarly for  $\beta$).  
 
 \begin{defn}(\textit{$\LL^p$}-ME)
 If,  for any $\lambda\in \Lambda$ and $\gamma\in \Gamma$, the integrals
$$\int_{X_\Lambda}\|\alpha(\gamma,x)\|_\Lambda^p d\mu(x)=\sum_{\lambda\in\Lambda} |\lambda|^p \mu(\{x:\alpha(\gamma,x)=\lambda\}) $$ 
$$\int_{X_\Gamma}\|\beta(\lambda,x')\|_\Gamma^pd\mu(x')=\sum_{\gamma\in\Gamma} |\gamma|^p \mu(\{x':\beta(\lambda,x')=\gamma\})$$
are finite, then the groups are called \textit{$\LL^p$-Measure Equivalent} ($\LL^p$-ME), where $\|\gamma\|_\Gamma 
$ denotes the distance of $\gamma$ from $\textbf{1}_\Gamma$ in the Cayley graph of $\Gamma$ (similarly, for $\Lambda$) and $1\leq p\leq \infty$. 
\end{defn}

 \vspace{2mm}

 \subsection{Mal'cev completion, Carnot Group, Carnot-Caratheodory Metric}\label{malcevcgcc}
 
 The following theorem is due to A. I. Mal'cev in \cite{Mal51}. 
 
  \begin{thm}\cite{Mal51}\label{Mal'cevcompletion}
 A group $\Gamma$ is isomorphic to a lattice in a simply connected nilpotent Lie group $N$ if and only if $\Gamma$ is finitely generated torsion free nilpotent group. 
\end{thm}
 
 \vspace{2mm}
 
\begin{defn}(\textit{Mal'cev completion}, \cite{Mal51})
The nilpotent group in the theorem is unique up to isomorphism with respect to the above-mentioned  countable Nilpotent group $\Gamma$. It is called the \textit{Mal'cev completion} of the nilpotent group $\Gamma$.  We denote it by $N(\Gamma)$. 
\end{defn}

We observe that 
$$\mbox{Rank}(\Gamma) \hspace*{1mm}=\hspace*{1mm} \mbox{Rank}\big(N(\Gamma)\big)\hspace*{1mm} =\hspace*{1mm} \mbox{dim}\big(N(\Gamma)\big) $$

 \vspace{2mm}

We consider any simply connected  Nilpotent Lie group $N$ and its Lie algebra  $\mathcal{N}$ and construct the \textit{Graded Lie Algebra} 

$$\mathcal{N}_\infty:=\oplus_{i=1}^r \big( \mathcal{N}^i/\mathcal{N}^{i+1}\big)$$ 

with the Lie bracket defined as follows: if $\bar{x}\in \mathcal{N}^{i}$ and $\bar{y}\in \mathcal{N}^{j}$

$$[\bar{x},\bar{y}]:=\bar{[x,y]},$$

which is well-defined modulo the subspace $\mathcal{N}^{i+j+1}$ (\cite{Pan83}). 

\vspace{2mm}
 
 \begin{defn}( \textit{Carnot group}, \cite{Pan83}) The simply connected (connected) graded nilpotent Lie group associated with the Lie algebra $\mathcal{N}_\infty$ is defined as the \textit{Carnot group} of the simply connected Lie group $N$. 

 \end{defn}
 
 \begin{defn}(\textit{Carnot group of a f.g. Nilpotent group}) If $\Gamma$ is a finitely generated (f.g.) torsion free nilpotent group, then the Carnot group associated to $\Gamma$ is defined as the Carnot group associated to its Mal'cev completion $N(\Gamma)$. We denote the Carnot group associated to $\Gamma$ by $\Gamma_\mathbb{R}$. 
 \end{defn}
 
  \vspace{2mm}
 
 \begin{defn}( \textit{Carnot-Caratheodory (CC) metric}, \cite{Pan83},\cite{BFT16})

Let $N$ be a nilpotent connected Lie group, and let $\mathcal{N}$ be its Lie algebra, and let $\mathcal{N}_{\mbox{sub}}$ be a vector subspace supplementing $[\mathcal{N},\mathcal{N}]$ equipped with a norm $\| \cdot \|$. A smooth path $\gamma : [0, 1] \rightarrow N$ is said to be \textit{horizontal} if $\big(\gamma(t)\big)^{-1}\gamma'(t)$ belongs to $\mathcal{N}_{\mbox{sub}}$ for all $t\in [0, 1]$. The length of $\gamma$ with respect to $\| \cdot \|$ is defined as 
 
 $$l(\gamma)=\int_{0}^{1} \|\big(\gamma(t)\big)^{-1}\gamma'(t)\| dt .$$ 
 
It can be shown that since $\mathcal{N}_{\mbox{sub}}$  generates the Lie algebra $\mathcal{N}$, every pair of points can be joined by a horizontal path. The \textit{Carnot-Caratheodory metric} associated to $\|\cdot\|$ is defined so that the distance between two points in $N$ is

$$d_{cc}(x, y) = \mbox{inf} \{l(\gamma) : \gamma(0) = x, \gamma(1) = y\},$$

where the infimum is taken over all piecewise horizontal paths $\gamma$. 
\end{defn}
 
 \begin{rem}(\textit{CC-Metric on a quotient of a Nilpotent Lie group}, \cite{BFT16})
 
 Suppose $p : N \rightarrow N'$ be a surjective homomorphism between nilpotent connected Lie groups, and let $d_{cc}$ be a Carnot-Caratheodory metric on $N$ associated to 
 $(\mathcal{N}_{\mbox{sub}}, \|\cdot\|)$,  Then, there exists a unique metric $d'_{cc}$ on $N'$ defined so that
 
$${d'}_{cc}(x', y') := \mbox{inf} \{d_{cc}(x, y) :  p(x) = x', p(y) = y'\} $$

This metric is easily seen to be the Carnot-Caratheodory Finsler metric on $N'$ associated to the normed vector space $(\mathcal{N'}_{\mbox{sub}}, \|\cdot\|')$ obtained by projecting $(\mathcal{N}_{\mbox{sub}}, \|\cdot\|)$.

\end{rem}

\begin{rem}\label{leftinvcc}
We remark that $d_{cc}$ (resp. $d'{cc}$) are both invariant under the left action of $N$ . 
\end{rem}

 \subsection{Equi-continuous and Distal topological system}\label{ecdistal}

Let $X$ be a compact Hausdorff topological space and $G$ be an infinite countable group acting on $X$ be homeomorphisms. We assume that $\mathcal{U}_{X}$ is the collection of open sets in $X$. 

\begin{defn} (\textit{Equi-continuous system})
The topological system $(X,G)$ is \textit{Equi-continuous}(EC) if for any neighbourhood $U\in\mathcal{U}_{X}$ of the diagonal $\Delta=\{(x,x): x\in X\} \subset X\times X$  $U$, there is neighbourhood $V\in\mathcal{U}_{X}$ of $\Delta$ such that whenever $(x_1,x_2)\in U$ , then $(g\cdot x_1, g\cdot x_2)\in V$ for all $g\in G$.

If $X$ is metrizable with metric $d$, then we can rephrase the definition as follows : 
for every $\epsilon > 0$, there exists $\delta > 0$ such that $d(x_1,x_2)<\delta$ implies $d(g x_1, g x_2)<\epsilon$ for all $g\in G $.
\end{defn}

\begin{defn}(\textit{Distal System})
The topological system $(X,G)$ is called \textit{Distal} if $x,y\in X$ with $x\neq y$, then there exists 
$U\in \mathcal{U}_X$ (depending on $x$ and $y$) such that $(g\cdot x, g\cdot y)\notin \mathcal{U}_X$
for all $g\in G$.

If $X$ is metrizable with metric $d$, then we can rephrase the definition as follows :  if 
$$ \inf_{g\in G} \hspace*{2mm}  d(g x, g x') > 0 $$ 
for all distinct $x,x'\in X$. 
\end{defn}

\begin{rem}
We remark that every distal  topological system is equi-continuous.
\end{rem}

\section{IME of nilpotent groups and coarse equivalence of box spaces}

In this section, we assume that $\Gamma$ and $\Lambda$ are two finitely generated countable torsion free nilpotent groups with symmetric generating set $S^\Gamma$ and $S^\Lambda$, respectively.

\subsection{IME of nilpotent groups and asymptotic cones}\label{asympcone}

M. Gromov proves in \cite{Gro93} that if two finitely generated groups are QI, then their \textit{asymptotic cones} are bi-Lipschitz homeomorphic.  In \cite{Aus16}, T. Austin generalizes this result for  nilpotent groups under the condition of being integrable measure equivalent ( see Theorem 1.1, p. 2, \cite{Aus16}).  Moreover, P. Pansu proves in \cite{Pan83} ( (10) Th\'{e}or\`{e}me Principal) and \cite{Pan89}  that the asymptotic cone of a nilpotent group $N$ is a simply connected nilpotent Lie group $N_{\mathbb{R}}$ , called `Carnot group' associated to $N$, equipped with Carnot-Caratheodory metric $d^N_{cc}$ (see Subsection \ref{malcevcgcc} for the definition of Carnot group and its Carnot-Caratheodory metric). Therefore, we obtain that for some $k'\geq 1$ there exist a $(k',1/k')$-bi-Lipschitz homeomorphism 

\begin{center}
\begin{tikzpicture}[node distance=3.5cm, auto] 

 \node (A) {$(\Gamma_{\mathbb{R}}, d^\Gamma_{cc})$};
 \node (B)[node distance=3.5cm, right of=A] {$(\Lambda_{\mathbb{R}},d^\Lambda_{cc})\hspace*{5mm}(*)$};
  \draw[->] (A) to node {$f$} (B);
 \end{tikzpicture}
\end{center}

\vspace{2mm}

\subsection{Scaling limit of finite quotients of nilpotent groups}

It is well known that a  nilpotent group is  \textit{residually finite}, i.e., there exists a decreasing sequence of finite index normal subgroups $\{N_n\}_{n\in\mathbb{N}}$ whose intersection is trivial ( $N=N_1\trianglerighteq N_2\trianglerighteq\cdots\trianglerighteq N_n\trianglerighteq\cdots\trianglerighteq{1}$). $\{N_n\}_{n\in\mathbb{N}}$  is called a \textit{residual chain} of $N$.  We consider a residual chain $\{\Gamma_n\}_{n\in\mathbb{N}}$ (resp. $\{\Lambda_n\}_{n\in\mathbb{N}}$ ) of $\Gamma$ (resp. $\Lambda$) satisfying the following properties: there exists an increasing sequence of positive integers $\{D^\Gamma_n\}_{n\in\mathbb{N}}$ (resp. $\{D^\Lambda_n\}_{n\in\mathbb{N}}$)  such that $B^\Gamma(1; D_n)$ is isometric to $B^{\Gamma/\Gamma_n}(1; D^\Gamma_n)$ (resp. $B^\Lambda(1; D^\Lambda_n)$ is isometric to $B^{\Lambda/\Lambda_n}(1;D^\Lambda_n)$) for all $n\in\mathbb{N}$.  Let $\Delta_n=\mbox{Max}\{\mbox{Diam}(\Gamma/\Gamma_n),\mbox{Diam}(\Lambda/\Lambda_n)\}$ for all $n\in\mathbb{N}$. The image of $S^\Gamma$ (resp. $S^\Lambda$) under the quotient map $\Gamma\twoheadrightarrow\Gamma/\Gamma_n$ (resp. $\Lambda\twoheadrightarrow\Lambda/\Lambda_n$) gives rise to a generating set  $S^\Gamma_n$ (resp. $S^\Lambda_n$) of $\Gamma/\Gamma_n$ (resp. $\Lambda/\Lambda_n$). 
We observe that the cardinalities of  $S^\Gamma_n$ and $S^\Lambda_n$ is uniformly bounded. 
Let $d^{\Gamma/\Gamma_n}$ (resp. $d^{\Lambda/\Lambda_n}$) be the Cayley graph distance in 
$\Gamma/\Gamma_n$ (resp. $\Lambda/\Lambda_n$) with respect to the generating set $S^\Gamma_n$ (resp. $S^\Lambda_n$). The following proposition says about the structure of the scaling limits of  the Cayley graphs of $\{(\Gamma/\Gamma_n, d^{\Gamma/\Gamma_n}/\Delta_n)\}$ and $\{(\Lambda/\Lambda_n, d^{\Lambda/\Lambda_n}/\Delta_n)\}$. 

\vspace{2mm}

We now state a proposition from \cite{BFT16} which gives a structure theorem of scaling limit of Cayley graphs of finite quotients of a Nilpotent group using a quantitative version of Gromov's theorem from \cite{BGT12}.

\begin{prop}(\cite{BFT16},p. 23, Theorem 3.5.1, Proposition 3.5.2)\label{scalinglimit}
 Suppose $\{N_n\}_{n\in\mathbb{N}}$ is a sequence of virtually nilpotent groups with bounded index and step, and bounded generating sets $T_n$. Then for every $\Delta_n\rightarrow\infty$, a subsequence of the Cayley graphs 
 $(N_n, d_{T_n} /\Delta_n, e)$ converges for the pointed Gromov-Hausdorff topology to a connected nilpotent Lie group, equipped with a left invariant polyhedral Carnot-Caratheodory metric.
\end{prop}

\vspace{2mm}
 
Therefore, after passing to a common subsequence, we obtain that  the Cayley graph of $\{(\Gamma/\Gamma_n, d_\Gamma/\Delta_n)$ (resp. the Cayley graph $\{(\Lambda/\Lambda_n, d_\Lambda/\Delta_n)$) converges for the pointed Gromov-Hausdorff Topology to a compact connected nilpotent Lie group, say $K_1$ (resp. $K_2$) with a left invariant Carnot-Caratheodory metric, say $d^{K_1}_{cc}$ (resp. $d^{K_2}_{cc}$). We refer the readers to Subsection \ref{malcevcgcc} for the definition of $d^{K_1}_{cc}$ (resp. $d^{K_2}_{cc}$).    Moreover,  by Lemma 2.2.3 in \cite{BFT16}, the Gromov-Haudorff Limit in both cases are finite dimensional torus, say $\mathbb{T}^{l_1}$ (resp. $\mathbb{T}^{l_2}$) with left invariant finsler metric (polyhedral Carnot-Caratheodory metric).

\subsection{IME between nilpotent groups and scaling limits}

By the proof Proposition of \ref{scalinglimit}, there exist projections $p_\Gamma: \Gamma_{\mathbb{R}}\rightarrow K_1$ (resp.  $p_\Lambda: \Lambda_{\mathbb{R}}\rightarrow K_2$) which is  $\Gamma_{\mathbb{R}}$-equivariant (resp. $\Lambda_{\mathbb{R}}$-equivariant ). Since  $\Gamma$ and $\Lambda$ are IME, using  $(*)$ in Subsection \ref{asympcone} we obtain a $(k',1/k')$-bi-Lipschitz homeomorphism 
\linebreak
$f: \Gamma_{\mathbb{R}}\rightarrow \Lambda_{\mathbb{R}}$ for some $k'\geq 1$. Using $f: \Gamma_{\mathbb{R}}\rightarrow \Lambda_{\mathbb{R}}$,  we can construct a $(k',1/k')$-bi-Lipschitz homeomorphism $\tilde{f}: K_1 \rightarrow K_2$ satisfying the following commutative diagram: 

 \begin{center}
\begin{tikzpicture}[node distance=2.5cm, auto]

  \node (A) {$\Gamma_{\mathbb{R}}$};
  \node (B) [below of=A] {$(K_1,d^{K_1}_{cc})\simeq \mathbb{T}^{l_1}$};
  \node (C) [node distance=4.5cm, right of=B] {$(K_2,d^{K_2}_{cc})\simeq \mathbb{T}^{l_2}$};
  \node (D) [node distance=4.5cm, right of=A] {$\Gamma_{\mathbb{R}}$};
  \draw[->] (A) to node {$f$} (D);
  \draw[->] (A) to node {$p_\Gamma$ (projection)} (B);
  \draw[->] (B) to node {$\tilde{f}$} (C);
  \draw[->] (D) to node {$p_\Lambda$ (projection)} (C);
\end{tikzpicture}
\end{center}

Therefore, we can assume that $K_1 \simeq \mathbb{T}^l$ and $K' \simeq \mathbb{T}^l$ for some $l\geq 1$. Hence, we have the following proposition 

\begin{prop}\label{bilipschitzscalinglimit}
If $\Gamma$ and $\Lambda$ are IME and the scaling limit of $\{(\Gamma/\Gamma_n, d_\Gamma/\Delta_n)\}$ and   $\{(\Lambda/\Lambda_n, d_\Lambda/\Delta_n)\}$  are  $(K_1, d^{K_1}_{cc})$ and $(K_2, d^{K_2}_{cc})$), respectively (which are homemorphic to same $\mathbb{T}^l$ for some $l\geq 1$ by the previous discussion), then there exists a $(k',1/k')$-bi-Lipschitz homeomorphism $\tilde{f}:(K_1, d^{K_1}_{cc}) \rightarrow (K_2, d^{K_2}_{cc})$. 
 \end{prop}
 
 \vspace{2mm}
 
 \subsection{Construction of coarse equivalence between box spaces:}

 A \textit{box space} of a finitely generated residually finite group $G$ with residual chain $\{G_n\}_{n\in\mathbb{N}}$ is defined as the disjoint union  $\sqcup_{n\in\mathbb{N}}G/G_n$ and denoted by $\square_{G_n}G$. We consider two  residually finite groups $\Gamma$ and $\Lambda$ with residual chains $\{\Gamma_n\}_{n\in\mathbb{N}}$ and $\{\Lambda_n\}_{n\in\mathbb{N}}$, respectively.  We say that there is a \textit{$(\rho_{+},\rho_{-},D)$-coarse equivalence} between $\square_{\Gamma_n}\Gamma$ and $\square_{\Lambda_n}\Lambda$ if there exists $\phi_n:\Gamma/\Gamma_n\rightarrow \Lambda/\Lambda_n$ such that 
 $$ \rho_{-}\big(d^{\Gamma/\Gamma_n}(x_1,x_2)\big) \leq d^{\Lambda/\Lambda_n}\big(\phi_n(x_1),\phi_n(x_2)\big)\leq  \rho_{+}\big(d^{\Gamma/\Gamma_n}(x_1,x_2)\big) $$
 for all $x_1, x_2\in \Gamma/\Gamma_n$ and the image of $\phi_n$ is $C$-dense in $\Lambda/\Lambda_n$, where $\rho_{+},\rho_{-}:[0,\infty)\rightarrow[0,\infty)$ are two monotonically increasing functions going to $\infty$ and $C\geq 0$ . 

\vspace{2mm}

By Theorem \ref{cecondn}, equivalently it suffices to prove that for all $n\in\mathbb{N}$ there exists 
$\phi_n:\Gamma/\Gamma_n\rightarrow\Lambda/\Lambda_n$ and 
$\psi_n:\Lambda/\Lambda_n \rightarrow \Gamma/\Gamma_n$ such that 

$$ d^{\Lambda/\Lambda_n}\big(\phi_n(x_1), \phi_n(x_2)\big)\leq \rho_{1}\big(d^{\Gamma/\Gamma_n}(x_1,x_2)\big),\hspace*{1mm}d^{\Gamma/\Gamma_n}\big(\psi(y_1), \psi(y_2)\big)\leq \rho_{2}\big(d^{\Lambda/\Lambda_n}(y_1,y_2)\big)$$

$$ d^{\Gamma/\Gamma_n}\big((\psi\circ\phi)(x),x\big)\leq D, \hspace*{5mm} d^{\Lambda/\Lambda_n}\big((\phi\circ\psi)(y),y\big)\leq D $$
 for all $x_1, x_2, x\in \Gamma/\Gamma_n$ and $y_1,y_2,y\in \Lambda/\Lambda_n$, where $\rho_1,\rho_2: [0,\infty) \rightarrow [0,\infty)$ are two monotonically increasing functions going to $\infty$ and independent of $n$ and $D$ is a positive constant.
  
\vspace{2mm}

 \begin{thm}\label{imece}
 If two countable nilpotent groups $\Gamma$ and $\Lambda$ are IME, then there exists a  coarse equivalence between the box spaces $\square_{\Gamma_n}\Gamma$ and $\square_{\Lambda_n}\Lambda$. 
  \end{thm}
  
  \begin{proof}
  
By Proposition \ref{bilipschitzscalinglimit},  $\{(\Gamma/\Gamma_n,d^\Gamma_n/\Delta_n)\}$ and $\{(\Lambda/\Lambda_n,d^\Lambda_n/\Delta_n )\}$ both converge to torus $\mathbb{T}^l$ for the pointed GH-topology as $n\rightarrow\infty$. Therefore,  for all $n\in\mathbb{N}$ we obtain $(1/n)$-isometries $q^\Gamma_n:\Gamma/\Gamma_n\rightarrow\mathbb{T}^l$ and  $q^\Lambda_n:\Lambda/\Lambda_n\rightarrow\mathbb{T}^l$. We fix $g\in \Gamma/\Gamma_n$ and consider $(\tilde{f}\circ q^\Gamma_n)(g)\in \mathbb{T}^l$.  Since the image of $q^\Lambda_n$ is $(1/n)$-dense in $ \mathbb{T}^l$ , there exist $h\in \Lambda/\Lambda_n$ such that $d\big((\tilde{f}\circ q^\Gamma_n)(g), q^\Lambda_n(h)\big)<1/n$. We define
 $$\phi_n(g):=h.$$
 Therefore, we obtain a map $\phi_n: \Gamma/\Gamma_n\rightarrow \Lambda/\Lambda_n$ such that the following diagram commute:

 \begin{center}
\begin{tikzpicture}[node distance=2.5cm, auto]

  \node (A) {$(\Gamma/\Gamma_n,d^\Gamma_n/\Delta_n )$};
  \node (B) [below of=A] {$(K_1,d^{K_1}_{cc})\simeq \mathbb{T}^{l}$};
  \node (C) [node distance=4.5cm, right of=B] {$(K_2,d^{K_2}_{cc})\simeq \mathbb{T}^{l}$};
  \node (D) [node distance=4.5cm, right of=A] {$(\Lambda/\Lambda_n,d^\Lambda_n/\Delta_n )$};
  \draw[->] (A) to node {$\phi_n$} (D);
  \draw[->] (A) to node {$q^\Gamma_n$ ($1/n$-isometry)} (B);
  \draw[->] (B) to node {$\tilde{f}$} (C);
  \draw[->] (D) to node {$q^\Lambda_n$ ($1/n$-isometry)} (C);
\end{tikzpicture}
\end{center}

We define $\rho_1:[0,\infty)\rightarrow[0,\infty)$ in the following way: $\rho|_{[0,D_1]}(r)=k'r+\Delta_1(k'+3)$ and $\rho|_{[D_i,D_{i+1}]}(r)=k'r+\Delta_{i+1}(k'+3)$ for all $i\in\mathbb{N}$. It follows from the definition that  $\rho_1$ is a monotonically increasing function going to $\infty$. 

\vspace{2mm}

\textit{Claim 1:}  $d^{\Lambda/\Lambda_n}\big(\phi_n(g),\phi_n(g')\big)\leq \rho_1\big(d^{
\Gamma/\Gamma_n}(g,g')\big)$ for all $g,g'\in \Gamma/\Gamma_n$.

\vspace{2mm}

\textit{Proof of Claim 1:}  We let $\phi_n(g)=h$ and $\phi_n(g')=h'$.   Since $q^\Lambda_n$ is a $(1/n)$-isometry, we have 
 
 $$| d^{K_2}_{cc}\big( q^\Lambda_n(h),q^\Gamma_n(h') \big) - d^{\Lambda/\Lambda_n}(h,h') | < 1/n$$
 which implies 
 
\begin{align*}
d^{\Lambda/\Lambda_n}(h,h') & \leq d^{K_2}_{cc}\big( q^\Lambda_n(h),q^\Gamma_n(h') \big) +\frac{1}{n}\\
 & \leq  d^{K_2}_{cc}\big( q^\Lambda_n(h),(\tilde{f}\circ q^\Gamma_n)(g)\big)+ d^{K_2}_{cc}\big( \tilde{f}\circ q^\Gamma_n)(g) , \tilde{f}\circ q^\Gamma_n)(g') \big)\\
 & d^{K_2}_{cc}\big(\tilde{f}\circ q^\Gamma_n)(g'), q^\Lambda_n(h')\big)+\frac{1}{n}
\end{align*}

By the definitions of $h$ and $h'$,  
$$d^{K_2}_{cc}\big( q^\Lambda_n(h),(\tilde{f}\circ q^\Gamma_n)(g)\big)<\frac{1}{n},\hspace*{2mm} d^{K_2}_{cc}\big(\tilde{f}\circ q^\Gamma_n)(g'), q^\Lambda_n(h')\big)
<\frac{1}{n}.$$
Since $\tilde{f}$ is  $(k',1/k')$-bi-Lipschitz, we obtain that 
$d^{K_2}_{cc}\big( \tilde{f}\circ q^\Gamma_n)(g) , \tilde{f}\circ q^\Gamma_n)(g') \big)\leq 
k' d^{K_1}_{cc}\big(q^\Gamma_n(g),q^\Gamma_n(g')\big)$. Therefore, we have 

$$d^{\Lambda/\Lambda_n}(h,h')\leq k' d^{K_1}_{cc}\big(q^\Gamma_n(g),q^\Gamma_n(g')\big)+\frac{3}{n}   $$

Again, since $q_n$ is a $(1/n)$-isometry, we have 

$$| d^{K_1}_{cc}\big( q^\Gamma_n(g),q^\Gamma_n(g') \big) - d^{\Gamma/\Gamma_n}(g,g') | < 1/n$$

which implies $d^{K_1}_{cc}\big( q^\Gamma_n(g),q^\Gamma_n(g') \big)< d^{\Gamma/\Gamma_n}(g,g')+\frac{1}{n}$. Therefore, we obtain that 

\begin{align*}
d^{\Lambda/\Lambda_n}(h,h') \leq k' d^{\Gamma/\Gamma_n}(g,g') +\frac{k'+3}{n}
                                                \leq \rho_1\big(d^{\Gamma/\Gamma_n}(g,g')\big)
\end{align*}

Hence, we prove our claim.  \hfill\(\Box\)

\vspace{2mm}

Let $\rho_2=\rho_1$. By a similar way, we can construct a map $\psi_n:\Lambda/\Lambda_n\rightarrow\Gamma/\Gamma_n$ such that 

$$d^{\Gamma/\Gamma_n}\big(\psi_n(h), \psi(h')\big)\leq \rho_{2}\big(d^{\Lambda/\Lambda_n}(h,h')\big)$$

for all $h,h'\in \Lambda/\Lambda_n$ and satisfying the following commutating diagram: 

 \begin{center}
\begin{tikzpicture}[node distance=2.5cm, auto]

  \node (A) {$(\Gamma/\Gamma_n,d^\Gamma_n/\Delta_n )$};
  \node (B) [below of=A] {$(K_1,d^{K_1}_{cc})\simeq \mathbb{T}^{l}$};
  \node (C) [node distance=4.5cm, right of=B] {$(K_2,d^{K_2}_{cc})\simeq \mathbb{T}^{l}$};
  \node (D) [node distance=4.5cm, right of=A] {$(\Lambda/\Lambda_n,d^\Lambda_n/\Delta_n )$};
  \draw[->] (D) to node {$\psi_n$} (A);
  \draw[->] (A) to node {$q^\Gamma_n$ ($1/n$-isometry)} (B);
  \draw[->] (C) to node {$(\tilde{f})^{-1}$} (B);
  \draw[->] (D) to node {$q^\Lambda_n$ ($1/n$-isometry)} (C);
\end{tikzpicture}
\end{center}

\vspace{2mm}

\textit{Claim 2:} There exists $D>0$ (independent of $n$) such that for all $n\in\mathbb{N}$,  $g\in\Gamma/\Gamma_n$ and $h\in\Lambda/\Lambda_n$

$$ d^{\Lambda/\Lambda_n}\big((\psi_n\circ\phi_n)(g),g\big)\leq D, \hspace*{5mm} d^{\Gamma/\Gamma_n}\big((\phi_n\circ\psi_n)(h),h\big)\leq D . $$

\vspace{2mm}

\textit{Proof of Claim 2:} Let $\phi_n(g)=g'$ and $\psi_n(g')=g''$. Therefore, we have 
$$d^{K_2}_{cc}\big(q^\Lambda_n(g'), (\tilde{f}\circ q^\Gamma_n)(g)\big)<\frac{1}{n}\hspace*{5mm}\dots (1),$$

$$d^{K_1}_{cc}\big(q^\Gamma_n(g''), (\tilde{f}^{-1}\circ q^\Lambda_n)(g')\big)<\frac{1}{n}\hspace*{5mm}\dots (2).$$

Since $d^{\Lambda/\Lambda_n}\big(q^\Lambda_n(g'), (\tilde{f}\circ q^\Gamma_n)(g'')\big)<\frac{1}{n}$ and $\tilde{f}$ is a $(k',1/k')$-bi-Lipschitz homeomorphism, using (1) and (2) we obtain that 

\begin{align*}
d^{\Gamma/\Gamma_n}\big(q^\Gamma_n(g), q^\Gamma_n(g'')\big)
  & \leq d^{\Gamma/\Gamma_n}\big(q^\Gamma_n(g), (\tilde{f}^{-1}\circ q^\Gamma_n)(g')\big)+ d^{\Gamma/\Gamma_n}\big((\tilde{f}^{-1}\circ q^\Gamma_n)(g'), q^\Gamma_n(g'')\big)\\
 & < \frac{k'}{n}+\frac{1}{n}
 \end{align*}
 
 Since the distortion of $q^\Gamma_n$ is $1/n$, we obtain that 
 
 $$ d^{\Gamma/\Gamma_n}(g,g'')< \frac{k'}{n}+\frac{1}{n} +\frac{1}{n} < \frac{k'+2}{n}<(k'+2). $$
 
 Let $D=k'+2$. Therefore, $d^{\Gamma/\Gamma_n}\big((\psi_n\circ\phi_n)(g),g\big)\leq D$ for all 
 $g\in \Gamma/\Gamma_n$ and for all $n\in\mathbb{N}$. Similarly, we obtain that 
 $d^{\Lambda/\Lambda_n}\big((\phi_n\circ\psi_n)(h),h\big)\leq D$ for all 
 $h\in \Lambda/\Lambda_n$ and for all $n\in\mathbb{N}$. 
Hence, we prove our claim. \hfill\(\Box\)

\vspace{2mm}

Therefore, $\square_{\Gamma_n}\Gamma$ and $\square_{\Lambda_n}\Lambda$ are coarsely equivalent.

\end{proof}

 \bigskip 
  
\begin{rem}\label{ballstoballs}
By Theorem  \ref{imece}, there exists a $(\rho_{+},\rho_{-},C)$-Coarse Equivalence between $\square_{\Gamma_n}\Gamma$ and $\square_{\Lambda_n}\Lambda$, $\rho_{+},\rho_{-}:[0,\infty)\rightarrow[0,\infty)$ are two monotonically increasing functions going to $\infty$ and $C\geq 0$. We denote it by the same notation $\{\phi_n:\Gamma/\Gamma_n\rightarrow\Lambda/\Lambda_n\}_{n\in\mathbb{N}}$. After passing to a subsequence, we can assume that $\phi_n\big(B^{\Gamma/\Gamma_n}(1;D^\Gamma_n)\big)\subset B^{\Lambda/\Lambda_n}(1;D^\Lambda_n)$ for all $n\in\mathbb{N}$ and $\phi_n(1)=1$ for all $n\in\mathbb{N}$. We observe that we can construct a $(\rho_{+},\rho_{-},C)$-Coarse Equivalence $\phi:\Gamma\rightarrow\Lambda$ using $\{\phi_n\}_{n=1}^\infty$ such that the following diagram commute: 

  \begin{center}
\begin{tikzpicture}[node distance=2.5cm, auto]

  \node (A) {$B^{\Gamma}(1;D^\Gamma_n)(\subset \Gamma$)};
  \node (B) [below of=A] {$B^{\Gamma/\Gamma_n}(1; D^\Gamma_n )
  (\subset \Gamma/\Gamma_n)$};
  \node (C) [node distance=5.5cm, right of=B] {$B^{\Lambda/\Lambda_n}
  (1; D^\Lambda_n)(\subset \Lambda/\Lambda_n)$};
  \node (D) [node distance=5.5cm, right of=A] {$B^\Lambda(1; D^\Lambda_n)
  (\subset \Lambda)$};
  \draw[->] (A) to node {$\phi$} (D);
  \draw[->] (A) to node {$\backsimeq$ (isometry)} (B);
  \draw[->] (B) to node {$\phi_{n}$} (C);
  \draw[->] (D) to node {$\backsimeq$ (isometry)} (C);
\end{tikzpicture}
\end{center}

\end{rem}

\section{Construction of distal topological coupling between nilpotent groups}

In this section make a simpler version of the construction of topological coupling done by the second author in \cite{Das18}. In Theorem \ref{topocoupling}, we construct a topological coupling $X$ between $\Gamma$ and $\Lambda$ such that $G=\Gamma\times\Lambda$ action on $X$ is distal. 

 \begin{thm}\label{topocoupling}
  If $\square_{\Gamma_n}\Gamma$ and $\square_{\Lambda_n}\Lambda$ are  coarsely equivalent, then there exists a topological coupling between $\Gamma$ and $\Lambda$, i.e., $\Gamma$ and $\Lambda$ have continuous  commuting faithful  co-compact properly discontinuous action on a locally compact  space $X$ ( equivalently $\Gamma$ and $\Lambda$ are coarsely equivalent).  Moreover,  $G=\Gamma\times\Lambda$ action on $X$ is distal and  minimal (hence, topologically transitive).  
\end{thm}. 

\vspace{2mm}

We will prove several lemmas which will lead to the proof of Theorem \ref{topocoupling}. Let  
$\square_{\Gamma_n}\Gamma$ and $\square_{\Lambda_n}\Lambda$ be coarsely 
equivalent, i.e., there exists a sequence of maps  $\phi_n:\Gamma/\Gamma_n\rightarrow \Lambda/\Lambda_n$  for all $n\in\mathbb{N}$ such that 
$$ \rho_{-}\big(d^{\Gamma/\Gamma_n}(x_1,x_2)\big) \leq d^{\Lambda/\Lambda_n}\big(\phi_n(x_1),\phi_n(x_2)\big)\leq  \rho_{+}\big(d^{\Gamma/\Gamma_n}(x_1,x_2)\big) $$
 for all $x_1, x_2\in \Gamma/\Gamma_n$ and the image of $\phi_n$ is $C$-dense in $\Lambda/\Lambda_n$, where $\rho_{+},\rho_{-}:[0,\infty)\rightarrow[0,\infty)$ are two monotonically increasing functions going to $\infty$ and $C\geq 0$.  By Remark \ref{ballstoballs}, we can assume that $\phi_n(D^\Gamma_n)\subset D^\Lambda_n $ and 
$\phi_n(1)=1$ for all $n\in\mathbb{N}$. In the set $\mathcal{F}_n$, defined as the collection of  all  $(\rho_{+},\rho_{-},C)$-Coarse Equivalence from $\Gamma/\Gamma_n$ to $\Lambda/\Lambda_n$ sending the identity to identity, we define $\Gamma$-action for all $n\in\mathbb{N}$ as follows: 

$$(\gamma\cdot \phi)([\gamma'])=  [\phi([\gamma]^{-1})]^{-1} \phi([\gamma]^{-1}[\gamma']).$$

For all $n\in\mathbb{N}$,  we consider the subset $Z_n$ defined as 

$$Z_n=\Gamma\cdot \phi_n \subset \mathcal{F}_n.$$

 We define the following metric $d_n$ on $Z_n$ for all $n\in\mathbb{N}$: 
 
$$d_n(\phi,\psi)=2^{-\mbox{sup}\{R | \phi|_{B_R(1_{\Gamma/\Gamma_n})}\equiv \psi|_{B_R(1_{\Gamma/\Gamma_n})}\}}.$$

\vspace{2mm}

\begin{lem}\label{relcompact}

 $\{(Z_n,d_n)\}_{n\in\mathbb{N}}$ is a sequence of compact metric spaces which are `relatively compact' in Gromov-Hausdorff Topology,i.e., 
\begin{itemize}
\item[(i)]  there is a constant $E$ such that $\mbox{diam} (Z_n) \leq E$ for all $n\in\mathbb{N}$;
\item[(ii)] for all $\epsilon>0$ there exists a natural number $N=N(\epsilon)$ such that every $X_n$ contains an $\epsilon$-net consisting of at most $N$ points. 
\end{itemize}

\end{lem}

\begin{proof}

 Fix $\epsilon >0$. We choose an arbitrarily large integer $R$ such that $2^{-R}<\epsilon$. Since $\mbox{diam}\hspace*{1mm}  Z_n \leq 1$ for all $n$, the first criterion is trivially satisfied. We now prove the second criterion. We construct a map 
$$h_n:Z_n\rightarrow \mathcal{F}\Big(B^{\Gamma/\Gamma_n}_R(1_{\Gamma/\Gamma_n}), B^{\Lambda/\Lambda_n}_{\rho_{1}(R)}(1)\Big)\hspace*{3mm} \mbox{ defined by}\hspace*{3mm} \xi \mapsto \xi|_{B^{\Gamma/\Gamma_n}_R(1)},$$
where $\mathcal{F}\Big(B^{\Gamma/\Gamma_n}_R(1), B^{\Lambda/\Lambda_n}_{\rho_{+}(R)}(1_{\Lambda/\Lambda_n})\Big)$ is the collection of all maps from the ball of radius $R$ in $\Gamma/\Gamma_n$ around $1_{\Gamma/\Gamma_n}$, denoted by $B^{\Gamma/\Gamma_n}_R(1)$,  to the ball of radius $\rho_{+}(R)$ in $\Lambda/\Lambda_n$ around 
$1$, denoted by $B^{\Lambda/\Lambda_n}_{\rho_{+}(R)}(1)$. Clearly, the map is well-defined. For all $Z_n$, we can  choose a subset $A_n$ of $Z_n$ such that the restriction of $h_n$ on $A_n$ is injective and the image of $A_n$ under $h_n$ covers the whole image of $Z_n$.

We claim that $A_n$ is an $\epsilon$-net in $X_n$.  Take any $\xi$ in $Z_n$. By the definition of $A_n$, there exists an element 
$\eta\in A_n$ so that $\xi$ and $\eta$ belong to the same fiber of the map $h_n$. Therefore, $\xi$ and $\eta$ coincides on the ball $B^{\Gamma/\Gamma_n}_R(1)$,
which implies that $d_{Z_n}(\xi,\eta)\leq 2^{-R}<\epsilon$. Hence, $A_n$ is an $\epsilon$-net in $Z_n$.

Now, we compute an estimate of the cardinality of the set $A_n$. We observe that 
\begin{eqnarray*}
\Big|\mathcal{F}\Big(B^{\Gamma/\Gamma_n}_R(1), B^{\Lambda/\Lambda_n}_{\rho_{+}(R)}(1)\Big)\Big| & \leq & | B^{\Gamma/\Gamma_n}_R(1)|\hspace*{1mm}|B^{\Lambda/\Lambda_n}_{\rho_{+}(R)}(1)|\\
                                                         & \leq & |S_G|^R\hspace*{1mm}|S_H|^{\rho_{+}(R)},
\end{eqnarray*}
where $S^\Gamma_n$ and $S^\Lambda_n$ are two generating subsets of $\Gamma/\Gamma_n$ and $\Lambda/\Lambda_n$, respectively. Therefore, $|A_n|\leq |S^\Gamma_n|^R\hspace*{1mm}|S^\Lambda_n|^{\rho_{1}(R)}$ for all $n$. Hence, $\{Z_n\}_{n\in\mathbb{N}}$  is uniformly totally bounded, which proves our lemma. 

\end{proof}

By the above-mentined lemma, after passing to a subsequence $\{(Z_n,d_n)\}_{n\in\mathbb{N}}$ converges to a compact metric space, say $(Z,d)$. We fix a non-principal ultrafilter $\omega$.  From \cite{BH99} ( Exercice 5.52), we obtain that  if a sequence of pointed metric spaces converges in the pointed GH topology to $X$, then its ultralimit with respect to $\omega$ is isometric to $X$,i.e., $(Z,d)$ is isometric to 
$\big(\Pi_n (Z_n,d_n,\phi_n)\big)/\sim$, where $(x_n)\sim (y_n)$ if $\lim_{n\rightarrow\omega} d_n(x_n,y_n)=0$.

\vspace{2mm}

 \begin{lem}\label{gpactiononZ}
There exist a continuous faithful $\Gamma$ action on $(Z,d)$.   
\end{lem}
\begin{proof}
Since  $(Z,d)$ is isometric to $(\Pi_n Z_n)/\sim$, where $\big((x_n)\big)\sim \big((y_n)\big)$ if 
\linebreak
$\lim_{n\rightarrow\omega} d_n(x_n,y_n)=0$, it suffices to define $\Gamma$-action on $(\Pi_n Z_n)/\sim$.  We fix $\gamma\in\Gamma$. We define $\gamma\cdot (x_n):=  (\gamma\cdot x_n)$. This action is well defined. Indeed, we let $\epsilon>0$. Without loss of generality, we assume that  $\epsilon< 1$. Let $\delta=\frac{\epsilon}{ 2^{|\gamma^{-1}|}}$. Suppose, $d_n(x_n,y_n)=2^{-R} <\delta$, where $R$ is a positive integer. Therefore, we obtain that
$$d_n(x_n,y_n)= d_n(\gamma^{-1}\gamma\cdot x_n,\gamma^{-1}\gamma\cdot y_n)=2^{-R}<\delta=\frac{\epsilon}{ 2^{|\gamma^{-1}|}}$$
which implies 
$d_n(\gamma\cdot x_n, \gamma \cdot y_n)\leq 2^{-(R-|\gamma^{-1}|)}<\epsilon$. 
 If $\big((x_n)\big)\sim \big((y_n)\big)$, $\omega(A_{\delta})=1$,  where $A_{\delta}=\{n: d_n(x_n,y_n)<\delta\}$. Let $B_{\epsilon}=\{n: d_n\big(\gamma\cdot x_n, \gamma\cdot y_n\big)< \epsilon  \}$. Since $A_\delta\subset B_\epsilon$, we get that $\omega(B_\epsilon)=1$. Therefore, $\Gamma$-action on $(Z,d)$ is well defined.  
 
 Now, we prove that the action is continuous. We fix $\gamma\in\Gamma$ and $x=\big((x_n)\big)\in Z$. Take $\epsilon<1$. Let $\delta=\frac{\epsilon}{2^{|\gamma^{-1}|}}$. If $d(x,y)=2^{-R}<\delta$
 for $y\in Z$, we obtain that 
 $$d(\gamma\cdot x,\gamma\cdot y)=2^{-(R-|\gamma^{-1}|)}<\epsilon$$. 
 Therefore, $\Gamma$-action on $(Z,d)$ is continuous.  
 
 Now, let $\gamma\neq 1$. Therefore, there exists $n_0\in\mathbb{N}$ such that $[\gamma]\neq 1$ in $\Gamma/\Gamma_n$ for all $n\geq n_0$. Without loss of generality, we can assume that there exists an infinite sequnce $\{n_k\}_{k=1}^\infty$ with $n_k\geq n_0$  such that  $\gamma$-action on $Z_{n_k}$ is not identity. Therefore,  there exist $\psi_{n_k}\in Z_{n_k}$  such that $\gamma\cdot \psi_{n_k}\neq \psi_{n_k}$. By Remark \ref{ballstoballs}, after passing through a subsequnce, we can construct $\psi\in Z$ such that $\gamma\cdot\psi\neq \psi$ (for detailed proof see Proposition 3.1  and Lemma 4.3 in \cite{Das18} ). Therefore, $\Gamma$ action on $(Z,d)$ is faithful. 
 
\end{proof}

\begin{lem}\label{equicont}
The actions of $\Gamma$ on $(Z_n,d_n)$ is equicontinuous for all $n\in\mathbb{N}$.  
\end{lem}

\begin{proof}

Let $E_n=\mbox{diam}(\Gamma/\Gamma_n, S^{\Gamma}_n)$.  We fix $\epsilon>0$. Without loss of generality, we assume that  $\epsilon<1$.
 Let  $\delta=\frac{\epsilon}{2^{E_n}}$.  If $d_n(\psi_n,\psi'_n)=2^{-R}<\delta$ for some positive integer $R$, 
$$d_n(\psi_n,\psi'_n)= d_n(\gamma^{-1}\gamma\cdot \psi_n,\gamma^{-1}\gamma\cdot \psi'_n)=2^{-R}<\delta=\frac{\epsilon}{2^{E_n}}$$
which implies 
$d_n(\gamma\cdot x_n, \gamma \cdot y_n)\leq 2^{-(R-|\gamma^{-1}|)}<\epsilon$  for all $\gamma\in\Gamma$. Since $\delta$ is independent  of $\gamma$, we obtain that $\Gamma$-actions of $(Z_n,d_n)$ is equicontinuous for all $n\in\mathbb{N}$.
\end{proof}

\vspace{2mm}

\begin{lem}\label{eqivproj}
There exists $\pi_n: Z\rightarrow Z_n$ for all $n\in\mathbb{N}$ such that for all $k\in\mathbb{N}$ there exists a subset $C_k$ of $\mathbb{N}$ with $\omega(C_k)=1$ such that 
\begin{itemize}
\item[(i)] $\pi_n$ is a $\frac{1}{k}$-isometry for all $n\in C_k$;
\item[(ii)]  $d_n\big(\pi_n(\gamma\cdot x),\gamma\cdot \pi_n(x)\big)<\frac{1}{k}$  for all $\gamma\in\Gamma$,  $x\in Z$ and $n\in C_k$. 
\end{itemize}
\end{lem}
\begin{proof}
We take $x\in Z= (\Pi_n Z_n)/\sim$. We consider a representative $(x_n)$ of $x$. 
We define $\pi_n(x)=x_n$. Let $\epsilon=\frac{1}{k}$. 
We take $x=[(x_n)]$ and $y=[(y_n)]$ in $Z$. Since $\lim_{n\rightarrow\omega} d_n(x_n,y_n)=d(x,y)$, we obtain that 
$$| d_n\big(\pi_n(x),\pi_n(y)\big) |=|d_n(x_n,y_n)-d(x,y)|<\frac{1}{k}$$
for all $n\in A_k$ with $\omega(A_k)=1$. 
On the other hand, if $[(x_n)]=[(x'_n)]$, then $\lim_{n\rightarrow\omega} d_n(x_n,x'_n)=0$, which implies there exists a subset $B_k$ of $\mathbb{N}$ with $\omega(B_k)=1$ such that for all $n\in B_k$
$$d_n(x_n,x'_n)<\frac{1}{k}. $$
Therefore, $\pi_n$ is a $\frac{1}{k}$-isometry for all $n\in C_k:=A_k\cap B_k$. 
Since $\gamma\cdot [(x_n)]=[(\gamma\cdot x_n)]$ for all $\gamma\in\Gamma$ and for all 
$[(x_n)]\in Z$, we have 
$$ d_n\big(\pi_n(\gamma\cdot (x_n)),\gamma\cdot \pi_n((x_n))\big)<\frac{1}{k}.$$
for all $n\in C_k$. Hence, we prove our lemma. 

\end{proof}

\begin{prop}\label{distalaction}
 $\Gamma$ action on $\Lambda\backslash X$ is distal. 
 \end{prop}

\begin{proof}
Let $\Gamma$ action $\Lambda\backslash X$ be not distal. Then, there exist $x=\big((x_n)\big)$
and $y= \big((y_n)\big)$  with $x\neq y$ and a sequence $\{\gamma_i\}_{i=1}^\infty$ in $\Gamma$ 
such that $\lim_{n\rightarrow\infty} \gamma_i \cdot x = \lim_{n\rightarrow\infty} \gamma_i \cdot y$.
Let $\epsilon=\frac{1}{k}$ for some positive integer grater than $k>1$. Then, there exists a positive integer $N_k$ such that  $d(\gamma_i\cdot x, \gamma_i\cdot y)<\frac{1}{k}$ for all $i\geq N_k$. Since 
$d(\gamma_i\cdot x, \gamma_i\cdot y)=\lim_{n\rightarrow\infty} d_n(\gamma_i\cdot x_n,  \gamma_i\cdot y_n)$, we obtain that there exists a subset $F_k$ of $\mathbb{N}$ with $\omega(F_k)=1$ such that 
$$ d_n(\gamma_i\cdot x_n,\gamma_i\cdot y_n)<\frac{1}{k}$$ 
for all $n\in F_k$ and for all $i\geq N_k$. Therefore, we have

\begin{align*}
d_n\big(\gamma_i\cdot \pi_n((x_n)), \gamma_i\cdot \pi_n((y_n))\big) & \leq d_n\big(\gamma_i\cdot \pi_n ((x_n)),  \pi_n (\gamma_i\cdot (x_n))\big)\\
& + d_n\big(\pi_n (\gamma_i\cdot (x_n)), \pi_n (\gamma_i\cdot (y_n))\big)+d_n\big(\pi_n (\gamma_i\cdot (y_n)), \gamma_i\cdot \pi_n((y_n))\big)\\
&<  \frac{2}{k}+d_n\big(\gamma_i\cdot (x_n)), \gamma_i\cdot (y_n)) \big)+\frac{1}{k} <\frac{4}{k}\hspace*{2mm}\dots (*)
\end{align*}
for all $i\geq N_k$ and for all $n\in H_k:=C_k\cap F_k$ with $\omega(H_k)=1$ ( using Lemma \ref{eqivproj}).  

Define $E_m=\mbox{Diam}(\Gamma/\Gamma_m)$ for all $m\in\mathbb{N}$. From the  inequality $(*)$, we obtain that 

$$d_n\big(\pi_n((x_n)),\pi_n((y_n))\big)<\frac{4}{k} 2^{|\gamma_i|}< \frac{4}{k} 2^{E_n}$$

We fix $\delta>0$. We  take sufficiently large $k$ such that  $\frac{4}{k} 2^{E_n}<\delta$. 
Therefore, there exists a subset $H_k$ of $\mathbb{N}$ with $\omega(H_k)=1$ such that for all $n\in H_k$

$$d_n\big(\pi_n((x_n)),\pi_n((y_n))\big)<\delta,   $$

which implies $\lim_{n\rightarrow\omega} d_n(x_n, y_n)=0$, in other words, 
$x=\big((x_n)\big)=\big((y_n)\big)=y$, which is  a contradiction. Hence, we prove our proposition.

\end{proof}

\begin{rem}
We remark that $\Gamma$ action on $(\Lambda\backslash X, d)$ can be also seen as 
a projective limit of the system of equicontinuous actions of $\Gamma$ on $(Z_n,d_n)$. We refer the reader to Definition \ref{projlimit} for projective limit.  Therefore, using Theorem 2.3 in \cite{Fur63}, we obtain that $\Gamma$ action on $(\Lambda\backslash X, d)$ is distal. 
\end{rem}

\begin{lem}\label{pttrans}
$\Gamma$-action on $(Z,d)$ is point transitive. 
\end{lem}
\begin{proof}

We consider  $\phi=\big((\phi_n)\big)$, where $\{\phi_n:\Gamma/\Gamma_n\rightarrow\Lambda/\Lambda_n\}_{n\in\mathbb{N}}$ is a collection of $(\rho_{+},\rho_{-},C)$-Coarse Equivalence
  ( $\rho_{+},\rho_{-}:[0,\infty)\rightarrow[0,\infty)$ are two monotonically increasing functions going to $\infty$ and $C\geq 0$ and they do not depend on $n$) and  $\phi_n\big(B^{\Gamma/\Gamma_n}(1;D^\Gamma_n)\big)\subset B^{\Lambda/\Lambda_n}(1;D^\Lambda_n)$ for all $n\in\mathbb{N}$ and $\phi_n(1)=1$ for all $n\in\mathbb{N}$ (using Theorem \ref{imece} and Remark \ref{ballstoballs}). We show that $\overline{\Gamma\cdot \phi}=Z$. Let $\psi=\big(\psi_n)\big)\in Z$. By definition of $Z_n$ $\phi_n=\gamma_n\cdot\phi_n$ for some $\gamma_n\in \Gamma$.   We consider the map $\pi_n:Z\rightarrow Z_n$. 
  Let $\epsilon=\frac{1}{k}$, where $k$ is an integer greater than 1,  As in  Lemma \ref{eqivproj}, We consider the subset $C_k$ of $\mathbb{N}$ with $\omega(C_k)=1$ and $\frac{1}{k}$-isometry 
$\pi_n:Z\rightarrow Z_n$  for all $n\in C_k$. Therefore, we obtain that 
$$|d(\psi,\gamma_n\cdot\phi)|<d_n(\psi_n,\gamma_n\cdot\phi_n)+\frac{1}{k} <\frac{1}{k}$$
for all $n\in C_k$ with $\omega(C_k)=1$. Therefore, we prove that $\overline{\Gamma\cdot \phi}=Z$, hence $\Gamma$-action on $(Z,d)$ is point transitive (hence topologically transitive). 
\end{proof}

\vspace{2mm}

\begin{rem}\label{ZisomZ'}
We consider the $(\rho_{+},\rho_{-},C)$-Coarse Equivalence $\phi:\Gamma\rightarrow\Lambda$ using  as obtained in Remark \ref{ballstoballs} such that the following diagram commute: 

  \begin{center}
\begin{tikzpicture}[node distance=2.5cm, auto]

  \node (A) {$B^{\Gamma}(1;D^\Gamma_n)(\subset \Gamma$)};
  \node (B) [below of=A] {$B^{\Gamma/\Gamma_n}(1; D^\Gamma_n )
  (\subset \Gamma/\Gamma_n)$};
  \node (C) [node distance=5.5cm, right of=B] {$B^{\Lambda/\Lambda_n}
  (1; D^\Lambda_n)(\subset \Lambda/\Lambda_n)$};
  \node (D) [node distance=5.5cm, right of=A] {$B^\Lambda(1; D^\Lambda_n)
  (\subset \Lambda)$};
  \draw[->] (A) to node {$\phi$} (D);
  \draw[->] (A) to node {$\backsimeq$ (isometry)} (B);
  \draw[->] (B) to node {$\phi_{n}$} (C);
  \draw[->] (D) to node {$\backsimeq$ (isometry)} (C);
\end{tikzpicture}
\end{center}

$Z'$ is the closure of the orbit  $\Gamma\cdot \phi$ inside the  set 

\begin{align*}
Z'' & =  \mbox{the collection of }\hspace*{2mm}(\rho_{+},\rho_{-},C)\mbox{-Coarse Equivalence} \hspace*{2mm}\psi:\Gamma\rightarrow\Lambda \\
 & \hspace*{2mm} \mbox{such that} \hspace*{2mm}\psi\big(B^{\Gamma}(1;D^\Gamma_n)\big)\subset B^\Lambda(1; D^\Lambda_n)\hspace*{2mm}\forall n\in\mathbb{N}\hspace*{2mm}\mbox{and}\hspace*{2mm} \psi(1)=1 
 \end{align*}

with the metric $d''$ defined as follows

$$d''(\psi_1,\psi_2)=2^{-\mbox{sup}\{R | \phi|_{B^\Gamma(1;R)}\equiv \psi|_{B^\Lambda(1;R)}\}}.$$

Clearly, $(Z'',d'')$ is a compact metrizable space ( the topology induced by $d''$ is same as point-wise convergent topology). Therefore, $Z'$ is a compact subset of $Z''$. We denote the induced metric on $Z'$ by $d'$. It is easy to see that that there exists a $\Gamma$-equivariant  isometry
$$ (Z,d) \longrightarrow (Z',d').$$
\end{rem}

\textit{Proof of Theorem \ref{topocoupling}:} 

We define the  topological space:
$$X=Z\times \Lambda,$$
with the metric  $d_X$ defined as follows: 
$$d_X\big((z_1,\lambda_1),(z_2,\lambda_2)\big)=d(z_1,z_2)+d_\Lambda(\lambda_1,\lambda_2)$$
for all $(z_1,\lambda_1), (z_2,\lambda_2) \in  Z\times \Lambda $, where $d$ is the metric on $Z$ and 
$d_\Lambda$ is the metric in the Cayley graph of $\Lambda$. 

It is easy to see that $X'':=Z''\times\Lambda$ equipped with the metric 
$$d_{X''}\big((z_1,\lambda_1),(z_2,\lambda_2)\big):=d''\big(z_1,z_2\big)+d_\Lambda(\lambda_1,\lambda_2)$$ 
where $(z_1,\lambda_1),(z_2,\lambda_2)\in X'$, is a closed subspace of 
$$ \Xi=\{f:\Gamma\rightarrow\Lambda : f \hspace*{2mm}\mbox{is a}\hspace*{2mm}(\rho_{+},\rho_{-},C)\mbox{-Coarse Equivalence}\}$$
with point-wise convergent topology and $\Gamma$ (resp. $\Lambda$) acts on this space by pre-multiplication (resp. post-multiplication) as follows: 
$$(\gamma\cdot f)(x)=f(\gamma^{-1} x)\hspace*{2mm}\mbox{and}\hspace*{2mm}(\lambda\cdot f)(x)=\lambda f(x)$$
for all $\gamma\in\Gamma$, $\lambda\in\Lambda$ and $x\in \Gamma$. Clearly,  $\Gamma$ and $\Lambda$ have commuting continuous actions on $X''$. Using Item $0.2. C_2'$ of \cite{Gro93}, we obtain that  $\Gamma$ and $\Lambda$ act  co-compactly, properly discontinuously on $X''$.

We define 
$$ X'=\overline{G\cdot\phi}\subset X''\subset \Xi .$$
 where $G=\Gamma\times\Lambda$. Using Remark \ref{ZisomZ'}, we extend the $\Gamma$-equivariant isometry 
$$ (Z',d')\longrightarrow (Z,d)$$ 
to an isometry  between $X'$ equipped with the restriction of the metric $d_{X''}$, denoted by $d_{X'}$,  and  $(X,d_X)$ by the following map $\Theta: X'\rightarrow X$ defined as follows: 
$$ \psi\mapsto \big(\tilde{\psi}, \psi(1)\big)$$
where $\tilde{\psi}(x)=\{\psi(1)\}^{-1}\psi(x)$ for all $x\in \Gamma$. 

From now onwards, we identify $ \psi$ with $\big(\tilde{\psi}, \psi(1)\big)$. We define $G=\Gamma\times\Lambda$ action on $X$ by the following way: 
$$(\gamma,\lambda)\cdot \big(\tilde{\psi}, \psi(1)\big):=\big(\gamma\cdot \tilde{\psi},\lambda \psi(1)).$$
It is easy to see that $\Theta$ is $G=\Gamma\times\Lambda$-equivariant. Since 
$\Gamma$ and $\Lambda$ act  co-compactly, properly discontinuously on the subspace $X'$ of $X''$, they act co-compactly, properly discontinuously on $X$ too. Therefore, $X$ is a topological coupling between $\Gamma$ and $\Lambda$. Proposition \ref{gpactiononZ} says that $\Gamma$ action on $(Z,d)$ is faithful, therefore $G$-action on $X$ is faithful too.

Now, we prove that $G$ action on $X$ is distal. We fix $\psi, \psi'\in X$ such that $\psi\neq \psi'$. 
  \begin{align*}
& \inf_{(\gamma,\lambda)\in G} d_X\big((\gamma,\lambda)\cdot (\tilde{\psi},\psi(1)), (\gamma,\lambda)\cdot (\tilde{\psi'},\psi'(1)) \\
& =  \inf_{(\gamma,\lambda)\in G} [d_Z\big(\gamma\cdot \tilde{\psi}, \gamma\cdot\tilde{\psi'}\big)+d_\Lambda\big(\lambda\psi(\gamma^{-1}), \lambda\psi'(\gamma^{-1}\big)]\\
& =  \inf_{\gamma\in \Gamma} d_Z\big(\gamma\cdot \tilde{\psi}, \gamma\cdot\tilde{\psi'}\big) + d_\Lambda\big(\psi(\gamma^{-1}), \psi'(\gamma^{-1}\big)\\
& \geq \inf_{\gamma\in \Gamma} d_Z\big(\gamma\cdot \tilde{\psi}, \gamma\cdot\tilde{\psi'}\big) \\
& > 0  \hspace*{5mm}\mbox{since}\hspace*{1mm}\Gamma\mbox{-action on Z is distal by Proposition \ref{distalaction}}
 \end{align*}

By Lemma \ref{pttrans}, $\Gamma$ action on $\Lambda\backslash X=Z$ is point transitive.  Now, using Lemma \ref{ptminimal} (Section 8), we obtain that  $\Gamma$ action on $\Lambda\backslash X$ is  minimal.
Therefore, $G=\Gamma\times\Lambda$-action on $X$ is distal and minimal (hence, topologically transitive).  

  \hfill\(\Box\)

\vspace{2mm}

\begin{rem}\label{quotientdistal}
 We remark that $\Gamma$-action on $\Lambda\backslash X$ and $\Lambda$-action on $\Gamma\backslash X$ are distal. $X$ is called a \textit{distal topological coupling} between $\Gamma$ and $\Lambda$ .  
\end{rem}

\vspace{2mm}

\section{Compactification of a coupling space and Furstenberg's theorem}\label{compacti_fur}

 Let $X$ be a topological coupling between the groups $\Gamma$ and $\Lambda$ as in Theorem \ref{topocoupling}. We denote by $\hat{X}=X\cup\{\infty\}$, the one point compactification of $X$.  Since $\hat{X}$ is a regular space with countable basis, using Urysohn Metrization Theorem, we can find a  metric, say $d_{\hat{X}}$, on $\hat{X}$ such that the restriction of $d_{\hat{X}}$ on $X$ gives the same topology as given by $d_X$. $G$ acts naturally by homeomorphisms on $(\hat{X}, d_{\hat{X}})$ fixing $\{\infty\}$. 
 
 \pagebreak

 \begin{lem}\label{distalcompacti}
 $G$ action on $(\hat{X}, d_{\hat{X}})$  is distal. 
 \end{lem}
 \begin{proof}
 Since $\overline{Gx}=X$ for all $x\in X$, we have 
 
$$ \inf_{g\in G}\hspace*{1mm} d_{\hat{X}}(gx,g\infty)=\inf_{g\in G}\hspace*{1mm} d_{\hat{X}}(gx,\infty) > 0 .$$
for all $x\in X$. Hence,  $G$ action on $(\hat{X}, d_{\hat{X}})$  is distal. 
 \end{proof}

\subsection{Furstenberg's Theorem}

\begin{defn}\label{projlimit}

A \textit{projective system of minimal flows} is a collection of minimal  flows (topological systems) $(Y_\lambda,H)$ indexed by ordinal numbers $\lambda\leq \theta$(for some ordinal $\theta$), and a
family of homomorphisms, $\phi_\nu^\rho: Y_\rho\rightarrow Y_\nu$, for $\nu<\rho\leq \theta$. 

\begin{itemize}
\item If $\nu < \rho < \eta \leq \theta$, then $\phi_{\nu}^\eta=\phi_{\nu}^\rho \circ \phi_{\rho}^\eta$. 
\item If $\mu\leq \theta $ is a limit ordinal , then $X_\mu$ is the minimal subset of the Cartesian flow $(\Pi_{\lambda<\mu}, H )$ consisting of all $y=(y_\lambda)\in \Pi_{\lambda<\mu} Y_\lambda $  with 
$y_\nu=\phi_\nu^\rho(y_\lambda)$ ($\nu<\rho<\mu$ ) and $\phi_{\rho}^\mu: Y_\mu\rightarrow Y_\lambda$ is the projection map. In this case, we say that $Y$ is the \textit{projective limit} of the flows $Y_\lambda$ ($\lambda<\mu$). 
\end{itemize}

\end{defn}

 \bigskip 
 
 \begin{defn}
 Let $(X,H)$ and $(Y,H)$ be two topological system and $\pi:X\rightarrow Y$ be a homo-morphism. Then, $\pi: X\rightarrow Y$ is called \textit{equi-continuous extension from $(Y,H)$ to $(X,H)$}  if for every $\alpha\in \mathcal{U}_X$ of the diagonal $\Delta_X$ in $X\times X$ there is a $\beta\in \mathcal{U}_X$ such that if $(x_1, x_2)\in \beta$ and $\pi(x_1)=\pi(x_2)$, then
 $(h x_1, h x_2)\in \alpha$ for all $h\in H$. 
 \end{defn}

The following structure theorem for distal  topological system for a compact space is due to H. Furstenberg ( see \cite{Fur63}, Page 483; \cite{Gla03} Theorem 1.31, p.34 )  

\begin{thm}\label{projlimitdistal}
Let $(Y,H)$ and $(Z,H)$ be distal topological systems and let $\pi:Y\rightarrow Z$ be a homomorphism. 
Then there is a projective system of minimal flows $(Y_\nu, G)$ ($\nu\leq \theta$) with $Y_\theta=Y $ and $Y_0=Z$ such that if $\nu<\theta$, $\phi_\nu^{\nu+1}: Y_{\nu+1}\rightarrow Y_\nu$ is an equi-continuous extension. In particular, if $(Y,H)$ is distal minimal, there is a projective  system as above such that $Y=Y_\theta$, the projective limit of this projective system, with $Y_0=1$ being the trivial system. 
\end{thm}

 \bigskip 

For a compact topological system $(Y,H)$ we define 
 $$Q_Y:=\cap_{\alpha\in \mathcal{U}_Y} \overline{H\alpha}$$
 where $\mathcal{U}_Y$ is the collection of entourage of the diagonal $\Delta_Y$ in $Y\times Y$. Clearly, $Q_Y$ is a $H$-invariant, closed, reflexive, symmetric relation on $Y$. Moreover, if $(Y,H)$  is minimal and distal,   $Q_Y$ is transitive (\cite{Gla03}, p.23). Hence, for a minimal distal compact topological system $(Y,H)$, $Q_Y$ is an equivalence relation on $Y$. 

\begin{thm}\label{maxecfactor}(\cite{Gla03} Theorem 1.10, p. 23)
 Let $(Y,H)$ be a distal minimal compact topological system as above. Then, $Y/Q_Y$  is the maximal  equi-continuous factor of $Y$. 
\end{thm}

\vspace{2mm}
 
 \subsection{Construction of Equicontinuous Coupling space}\label{ec-coupling}

We consider the set-up of Theorem \ref{topocoupling}. Let $\hat{X}=X\cup\{\infty\}$ be the one point compactification of $X$ with the metric $d_{\hat{X}}$ as obtained in 
Section \ref{compacti_fur}.  By Section \ref{compacti_fur}, $ \hat{X}/Q_{\hat{X}}$ is  the maximal equi-contnuous  minimal topological system  of the distal minimal topological system $(\hat{X}, G)$. We define the metric $\bar{d}$ on $ \hat{X}/Q_{\hat{X}}$ as follows which provides the quotient topology on $\hat{X}$ : 

$$ \bar{d}([x],[y]):=\mbox{inf}\{d_{\hat{X}}(\hat{x},\hat{y}) : \hat{x},\hat{y}\in \hat{X}\hspace*{1mm}\mbox{and}\hspace*{1mm} [x]=[\hat{x}], [y]=[\hat{y}]\}\dots (1) $$\label{metricdefn}
for all $[x],[y]\in \hat{X}/Q_{\hat{X}}$ . 

\vspace{2mm}

Since $G$-action on $\big(\hat{X}/Q_{\hat{X}},\bar{d}\big)$ is equi-continuous, we define a metric $\bar{D}$ on $\hat{X}/Q_{\hat{X}}$ 
as follows: 
$$\bar{D}([x_1],[x_2])=\hspace*{1mm}\inf_{g\in G} \bar{d}(g [x_1], g [x_2])$$
for all $[x_1],[x_2]\in \hat{X}/Q_{\hat{X}}$. We observe that $\bar{D}$ is an equivalent $G$-invariant metric on $\hat{X}/Q_{\hat{X}}$ (see Exercise 1.4, p. 18, \cite{Gla03}).

\vspace{2mm}

We observe that  the class of $\{\infty\}$ in $Q_{\hat{X}}$ consists of elements from the set 
$X\times\{\infty\}\cup \{\infty\}\times X \cup \{(\infty,\infty)\}$ and there is a $G$-equivariant homeomorphism  
$$\hat{X}/Q_{\hat{X}}-\{[\infty]\}\hspace*{1mm}\longrightarrow \hspace*{1mm}  X/Q_X$$
where $Q_X=\cap_{\beta\in \mathcal{U}_X} \overline{ G\beta} $, where $\mathcal{U}_X$ is the collection of entourage of the diagonal $\Delta_X$ in $X\times X$. We define 
$$ X_\II:= X/Q_X $$

\vspace{2mm}

 We induce the restriction of the metric $\bar{D}$ on $\hat{X}/Q_{\hat{X}}-\{[\infty]\}$ to $X_\II$ by the above-mentioned homeomorphism and denote it by $d_{X_\II}$. Therefore, we obtain the following proposition:
 
 \begin{prop}\label{isomcoup}
 There  is a topological coupling $\big(X_\II, d_{X_\II})$ between $\Gamma$ and $\Lambda$ such that $G=\Gamma\times\Lambda$ acts  on $X_\II$  faithfully, minimally (hence, topologically transitively) and isometrically (hence, equi-continuously). 
 \end{prop}
 
 \vspace{2mm}

\begin{rem}\label{quotientisom}
 $\Lambda$ action on $\Gamma\backslash X_\II$ and $\Gamma$ action on $\Lambda\backslash X_\II$ are isometric (hence, equi-continuous). $X_\II$ is called an \textit{isometric topological coupling} ( or \textit{equi-continuous topological coupling}) between $\Gamma$ and $\Lambda$. 
\end{rem}

\vspace{2mm}

\begin{rem}\label{propermap}
We remark that the quotient map $X\rightarrow X/Q_X= X_\II$ is a proper map, since it is the restriction map of a continuous function between one-point compactification 
of $X$ and $X_\II$. 
\end{rem}

\vspace{2mm}

\section{R. Ellis' enveloping monoid and EC-coupling space}

Let  $(Z,H)$ be a compact topological system. We consider the closure of the image of $H$ under the map 
$H \subset \mbox{Homeo}(Z) \subset Z^Z$, where the latter space is taken with the product topology, that is the topology of pointwise convergence. This is called the \textit{enveloping monoid} associated to the $H$-action on  $Z$. We denote it by $E(Z,H)$. Now, we state a theorem providing several equivalent conditions of a compact topological system being distal (resp. equi-continuous) in terms of its enveloping monoid.

\begin{thm}\label{envmonoiddistal}(\cite{Glas93}, Theorem 1.7, page 22)
 The following conditions are equivalent: 
 \begin{itemize}
 \item[1.] $(Z, H)$ is distal; 
 \item[2. ] $E(Z, H)$ is a group of homeomorphisms of $Z$. 
 \end{itemize}

 \vspace{2mm}

\end{thm} 
\begin{thm}\label{envmonoidequi}(\cite{Glas93}, Theorem 1.8, page 22)
  The following conditions are equivalent:
\begin{itemize}
\item[1.] $(Z, H)$ is EC;
\item[2.] $E(Z, H)$ is a group of homeomorphisms of $Z$ and the topologies of pointwise and uniform convergence coincide on $E(Z, H)$;
\item[3.] $E(Z, H)$ is a topological group whose action on $Z$ is jointly continuous; 
\item[4.] $E(Z,H)$ is a group of homeomorphisms of $Z$. 
\end{itemize}
 \end{thm}
 
 \vspace{2mm}
 
 \begin{lem}\label{ptminimal}
If a compact topological system $(Y,H)$ is distal and point transitive, i.e., there exists $y_0\in Y$ such that the orbit $Hy_0$ is dense in $Y$, then it is minimal, i.e., $Hy$ is dense in $Y$ for all $y\in Y$.  
\end{lem}
\begin{proof}
By Exercise 1.26(1) of \cite{Gla03}, $Y=\overline{Hy_0}=Ey_{0}$, where $E$ is Ellis' enveloping monoid associated to $(Y,H)$.  By Exercise 1.26 (2) of \cite{Gla03}, $Y=\overline{Hy_0}$ is minimal if and only if for any minimal ideal $M$ in $E$ ( see Exercise 1.23 \cite{Gla03} for the definition of \textit{minimal ideal}) $\overline{Hy_0}=My_0$. Since $(Y,H)$ is distal, $E$ is a group ( see Theorem \ref{envmonoiddistal}),  which implies any minimal ideal $M=E$. Therefore, $\overline{Hy_0}=Ey_0=My_0$. Hence, $(Y,H)$ is minimal. 

\end{proof}
 
\vspace{2mm}

\subsection{Elli's enveloping monoid as Nilpotent locally compact group }\label{envmonoid}

   Let  $\widehat{X_\II}=X_\II\cup\{\infty\}$, where $X_\II$ is equipped with the metric $d_{X_\II}$.  We consider the $G$-action on  $\widehat{X_\II} \simeq \hat{X}/Q_{\hat{X}}$ equipped with the metric $\bar{D}$. as in Subsection \ref{ec-coupling}. Let $\hat{E}$ be the enveloping monoid associated to the $G$-action on  $\widehat{X_\II}$.   By Theorem \ref{envmonoidequi},  $\hat{E}$ is a  topological group with uniform convergence topology. By the usual abuse of notation, we denote by 
  $\infty \in  \widehat{X_\II}^{\widehat{X_\II}}$ the corresponding constant map and we set $E = \hat{E} - \{\infty\}$.  It is easy to see that $E$ is  a locally compact space containing the countable  nilpotent group $G$ as a dense subset.

\vspace{2mm}

Now, we state a Proposition \ref{ecactellisgp} which is a generalization of  a result on abelian groups in \cite{Ba} to nilpotent groups. Since any two elements commute in an abelian group, the statement is relatively easier in \cite{Ba}. 

\begin{prop} \label{ecactellisgp} 
Let $(Y,d_Y)$ be a locally compact second countable  metric space, $\Gamma$ and $\Lambda$ be two finitely generated torsion free nilpotent group acting faithfully, isometrically, co-compactly and properly discontinuously on $Y$ and the action of $H=\Gamma\times\Lambda$ on $Y$ be minimal. Then, there exists a locally compact nilpotent group $F$, a continuous injective homomorphism $H\hookrightarrow F$ and a compact subgroup $C$ satisfying $\Gamma\cap C=\Lambda\cap C=\{1\}$ and $H$-equivariant homeomorphism $F/C\rightarrow Y$. 
\end{prop}

\begin{proof}

Let $\hat{Y}=Y\cup\{\infty\}$ be the one-point compactification of $(Y,d_Y)$. 
We consider Ellis' enveloping monoid $E(\hat{Y}, H)$. By Theorem \ref{envmonoidequi}, $E(\hat{Y}, H)$ is a  topological subgroup of $\mbox{Homeo}(\hat{Y})$ equipped with point-wise convergent topology (equivalently, with compact-open topology). We define 
$$F= E(\hat{Y}, H) - \{\infty\}$$ 
where $\infty\in \hat{Y}^{\hat{Y}}$ denotes the constant map to $\infty$. 

\vspace{2mm}

Since $H$ acts minimally on $Y$, we have 
$$ \overline{Hy}= Y$$
for all $y\in Y$. Therefore, by Exercise 1.26 (1) in \cite{Gla03} (p. 31) 
$$Fy=\overline{Hy}= Y,$$
which implies that for  every $p \in F$, $p(Y) \subset Y$ and the action of $F$ on $Y$ is transitive. Hence, $F$ is a collection of functions $Y\rightarrow Y$ such that for $p,q\in F, pq\neq\infty$. By Theorem \ref{envmonoidequi},  the product map $F\times F \rightarrow F$ and the action map $F \times Y \rightarrow Y$  are continuous. Since $H$ acts isometrically on $(Y, d_Y)$ and $H$ is a dense subgroup of $F$, by Theorem \ref{envmonoidequi}, $F$ is a locally compact nilpotent subgroup of the topological group $\mbox{Isom}\big((Y,d_Y)\big)\subset \mbox{Homeo}(Y)$ equipped with the compact-open topology. 

\vspace{2mm}

Now, if $F$ acts freely on $Y$, we will take $C=\{1\}$ and we define $\sigma: F\rightarrow Y$ by 
$$p\mapsto py$$
 for any fixed point $y\in Y$.  Otherwise, we fix a point $y_0\in Y$ and consider 
the subgroup $C:=Fix(y_0)$ of $F$, defined as the set of elements which fixes $y_0$,  containing a non-identity element. We consider the following compact subgroup of $F$
$$C'= [\Pi_{h\in H} \hspace*{2mm}\mbox{Isom}(\overline{B(y_0,R_h)},y_0,d_Y) ]\cap F,$$
where $\mbox{Isom}(\overline{B(y_0,R_h)},y_0,d_Y)$ is the isometry group of the closed ball $\overline{B(y_0,R_h)}$ around $y_0$ and fixing $y_0$ and $R_h=d_Y(y_0,hy_0)$ for all $h\in H$.   Since $C$ is a closed subgroup of the compact group $C'$. Therefore, $C$ is  a compact subgroup of $F$. 
We define $\sigma: F/C \rightarrow Y$ as follows: 
$$[p]\mapsto py_0 $$ for all $p\in F$.

\vspace{2mm}

Therefore, $\sigma$ is an $H$-equivariant continuous bijection. The point-wise convergence topology on $F$ is generated by the $F$-translations of sets of the form $\{p \in F : p y_0 \in O\}$ where $O$ is an open subset of
 $Y$. As the image of such a set under the orbit map is $O$, it is open in $Y$. We conclude that $\sigma$ is a $H$-equivariant homeomorphism. By definition of $F$, $H$ is a dense subgroup of $F$. Since $\Gamma$ (resp. $\Lambda$) acts on properly discontinuously on $Y$, hence on $F/C$, it will act properly discontinuously on $F$ too. Therefore, $\Gamma$ (resp. $\Lambda$) is a discrete subgroup of $F$. Since 
 $\Gamma$ (resp. $\Lambda$) is a discrete torsion free group, we have 
 $$\Gamma\cap C=\Lambda\cap C=\{1\}.$$
 Since $\Gamma$ (resp. $\Lambda$) acts co-compactly on $Y$ , hence on $F/C$, it will act on $F$ co-compactly too (as $C$ is a compact space).   Hence, we prove our proposition.

\end{proof}

\vspace{2mm}

We put $Y=X_\II$, $F=E$ and $H=G$ in Proposition \ref{ecactellisgp}. Since $(X_\II, d_{X_\II})$ is a locally compact, second countable metric space,  using Proposition \ref{ecactellisgp} we obtain the following corollary: 

\begin{cor}\label{eccouplingnp}
We consider the isometric coupling $(X_\II, d_{X_\II})$ of the finitely generated nilpotent  groups $\Gamma$ and $\Lambda$. Then, there exists a locally compact nilpotent group $E \subset \mbox{Isom}(X_\II)$, a continuous injective homomorphism $G=\Gamma\times\Lambda \hookrightarrow E$ with $\Gamma$ and $\Lambda$ being co-compact lattices in $E$ and 
a compact subgroup $C$ satisfying $\Gamma\cap C=\Lambda\cap C=\{1\}$ and a $G$-equivariant homeomorphism $E/C\rightarrow X_\II$, i.e., there exists an isomorphism of coupling spaces
$$  E/C \simeq  X_\II.$$ 
\end{cor}

\section{Proof of our main theorem}

The following is a deep structure theorem for locally compact, Hausdorff, connected topological group. 
 
 \begin{thm}(Gleason-Yamabe)\cite{Yam53}\label{gleyam}
 Let $H$ be locally compact, Hausdorff, connected topological group. Then there is a normal compact subgroup $K$ of $H$ such that the quotient $H/K$ is a Lie group.
\end{thm}

 \vspace{2mm}
 
 Using this theorem, we obtain the following proposition: 
 
 \begin{prop}\label{nilliegp}
There exists a connected nilpotent Lie group $E'$ containing  $\Gamma$ and $\Lambda$ as co-compact lattices and $G=\Gamma\times\Lambda$ as a nilpotent countable dense subgroup.  
 \end{prop}
 
 \begin{proof}
 We consider the connected component of the identity of the locally compact, Hausdorff group $E$ as obtained in Proposition \ref{ecactellisgp}. Since the connected component of a topological group  ( with Hausdorff condition ) is a closed subgroup of the topological group and $E$ contains $G$ as a dense subgroup, we obtain that $E$ is connected. Now, by Gleason-Yamabe Theorem \ref{gleyam}, there is a normal compact subgroup $K$ of $E$ such that the quotient $E'=E/K$ is a Lie group.  $E$ is a connected nilpotent group implies that $E'$ is a connected nilpotent Lie group too. Since $\Gamma$ and $\Lambda$ are discrete and have no torsion elements, we have $\Gamma\cap E=\{1\}$ and $\Lambda\cap E=\{1\}$. Therefore, the quotient map $E\rightarrow E'$ is injective on $\Gamma$ and $\Lambda$, respectively. We identify  $\Gamma$ (resp. $\Lambda$) in $E$ with its image in $E'$ under the quotient map $E\rightarrow E'$.   Hence, we obtain a nilpotent Lie group $E'$ such that it contains  $\Gamma$ and $\Lambda$ as co-compact lattices and $G=\Gamma\times\Lambda$ as a nilpotent countable dense subgroup.  
 \end{proof}
 
 \begin{rem}\label{maximaltorus}
 Let $\mathbb{T}^d$ be a maximal (central) torus in $E'$. Therefore, the compact subgroup  $C':=C/K\simeq C/(K\cap C)$ of $E'$ is a central torus inside $\mathbb{T}^d$.  Since $\Gamma$ and $\Lambda$ are discrete subgroup of $E'$ without torsion elements, we have 
 $$ \Gamma\cap C'=\Lambda\cap C'=\{1\}.$$
 \end{rem}

  \vspace{2mm}
  
  Now, we prove our main theorem \ref{MainTheorem}. 
  
\vspace{2mm}

 \textbf{Proof of Theorem \ref{MainTheorem}:} 

\vspace{2mm}

Let $\Gamma$ and $\Lambda$ be two finitely generated torsion free nilpotent groups which are IME. Using Theorem \ref{imece} and Theorem \ref{topocoupling}, we obtain that there exists a locally compact space such that
 
\begin{itemize}
\item[(i)] $\Gamma$ and $\Lambda$ act continuously, faithfully, co-compactly and properly discontinuously on $X$;
\item[(ii)] $\Gamma$ and $\Lambda$ actions on $X$ commute; 
\item[(iii)] $G=\Gamma\times\Lambda$-action on $X$ is minimal ( hence, topologically transitive) ;
 \item[(iv)] $G=\Gamma\times\Lambda$-action on $X$ is distal and the induced action of  $\Gamma$ (resp. $\Lambda$ ) on $\Lambda\backslash X$ (resp. on $\Gamma\backslash X$ ) is distal too. 
\end{itemize}

 Let $\hat{X}=X\cup\{\infty\}$. Since the induced action of $G$ on $\hat{X}$ is minimal and distal,  Theorem \ref{maxecfactor} provides the maximal equi-continuous factor  $\pi: \hat{X} \rightarrow \hat{X}/ Q_{\hat{X}}$ of  $\hat{X}$. We consider  $X_\II=X/Q_X  \simeq \hat{X}/ Q_{\hat{X}}-\{[\infty]\}$, where $Q_X$ is a $G$-invariant closed equivalence relation on $X$, with the metric $d_{X_\II}$ as described in Subsection 5.2.  By Proposition \ref{isomcoup},  $X_\II$  is a locally compact space satisfying the following properties: 

\begin{itemize}
\item[(a)] It is easy to see that $\Gamma$ and $\Lambda$ act continuously, co-compactly and properly discontinuously (since the quotient map from $X$ to $X_\II$ is proper by Remark \ref{propermap})  on $X_\II$. $G=\Gamma\times\Lambda$ acts faithfully on $X_\II$ since the map $\pi: \widehat{X}\rightarrow \widehat{X_\II}$, where $\widehat{X}=X\cup\{\infty\}$ and  $\widehat{X_\II}=X_\II\cup\{\infty\}$,
is equi-continuous extension of faithful, equi-continuous action of $G$ on $\widehat{X_\II}$;
\item[(b)] $\Gamma$ and $\Lambda$ actions on $X_\II$ commute; 
\item[(c)] $G=\Gamma\times\Lambda$ acts on $X_\II$  minimally  (hence, topologically transitively);
\item[(d)] $G=\Gamma\times\Lambda$-action on $X_\II$ is isometric and the induced action of  $\Gamma$ (resp. $\Lambda$ ) on $\Lambda\backslash X$ (resp. on $\Gamma\backslash X$ ) is isometric too. 
\end{itemize}

By  Corollary \ref{eccouplingnp}, $X_\II\simeq E/C$ as a coupling space between $\Gamma$ and $\Lambda$, where $E$ is a locally compact nilpotent group, $C$ is a compact subgroup of $E$ satisfying $\Gamma\cap C=\Lambda\cap C=\{1\}$ and there exist injective homomorphisms from $\Gamma$ and $\Lambda$ to $E$ such that the images of $\Gamma$ and $\Lambda$ inside  $E$ are co-compact lattices in $E$.

\vspace{2mm}

Now, using Proposition \ref{nilliegp}, we obtain  a connected nilpotent Lie group $E'$  with $\mathbb{T}^d$  being the maximal compact torus ($d\geq 0$ is the dimension of the torus) and containing  $\Gamma$ and $\Lambda$ as co-compact lattices and $G=\Gamma\times\Lambda$ as a  countable nilpotent dense subgroup.  We consider a connected, simply connected covering $\tilde{E'}$ of the connected Nilpotent Lie group $E'$ and let $f:\tilde{E'}\rightarrow E'$ be the covering map.  Since $\Gamma$ (resp. $\Lambda$) is discrete and has no torsion elements, we have $\Gamma\cap \mathbb{T}^d=\{1\}$ (resp. $\Lambda\cap \mathbb{T}^d =\{1\})$ and the compact subgroup $C'=C/K\simeq C/(K\cap C)\subset \mathbb{T}^d$. Therefore, we have the following central extensions: 

$$ 1\rightarrow \mathbb{Z}^d\rightarrow \tilde{E'} \rightarrow E' \rightarrow 1$$

$$ 1\rightarrow \mathbb{Z}^d\rightarrow \tilde{\Gamma} \rightarrow \Gamma \rightarrow 1$$

$$ 1\rightarrow \mathbb{Z}^d\rightarrow \tilde{\Lambda} \rightarrow \Lambda \rightarrow 1$$
where $\tilde{\Gamma}$ (resp.  $\tilde{\Lambda}$) are pre-images of $\Gamma$ (resp. $\Lambda$) under the map $f$. We observe that $\tilde{\Gamma}$ (resp.  $\tilde{\Lambda}$) are co-compact lattices in $\tilde{E'}$.
 
Since $\tilde{E'}$ and $N(\Gamma)$ (resp. $N(\Lambda)$) are simply connected and  $\tilde{\Gamma}$ (resp. $\tilde{\Lambda}$) is a co-compact lattice in $\tilde{E'}$, the restriction  map of  $f$ from $\tilde{\Gamma}$ to  $\Gamma$ ( resp. $\tilde{\Lambda}$ to $\Lambda$)   extends to a homomorphism  $\chi_\Gamma: \tilde{E'} \rightarrow N(\Gamma)$  (resp. $\chi_\Lambda: \tilde{E'} \rightarrow N(\Lambda)$) ( see \cite{Mal51}). Moreover, since there are no proper closed, connected subgroups of a simply connected nilpotent Lie group containing a co-compact lattice (\cite{Mal51}),  we obtain that $\chi_\Gamma$ (resp. $\chi_\Lambda$ ) is surjective.

\begin{center}
\begin{tikzpicture}[node distance=2cm, auto] 
 \node (A) {$N(\Gamma)$};
 \node (B)[node distance=3cm, right of=A] {$\tilde{E'}$};
 \node (C)[node distance=3cm, right of=B]{$N(\Lambda)$};
 \node (D)[below of=B]{$E'$};
   \draw[->] (B) to node {$\mathcal{\chi}_\Gamma $} (A);
 \draw[->] (B) to node {$\mathcal{\chi}_\Lambda $} (C);
 \draw[->] (B) to node {$f $} (D);

 \end{tikzpicture}
\end{center}

 Let $K_\Gamma=\mbox{Ker} (\chi_\Gamma)$ and $K_\Lambda=\mbox{Ker} (\chi_\Lambda$). $K_\Gamma$ (resp. $K_\Lambda$) is a closed Lie subgroup of $\tilde{E'}$.  We note that $\mathbb{Z}^d$ is a central subgroup of $K_\Gamma$ (resp. of $K_\Lambda$).  Since $\Gamma$ and $\Lambda$ are IME, by Remark \ref{qigrowth} we have 
 $$ \mbox{dim}\big(N(\Gamma)\big)=\mbox{Rank}\big(N(\Gamma)\big)=\mbox{Rank}(\Gamma)=\mbox{Rank}(\Lambda)=\mbox{Rank}\big(N(\Gamma)\big)=\mbox{dim}(N(\Gamma)).$$
  Therefore, we obtain that 
 $$\mbox{dim}\big( K_\Gamma)=\mbox{dim}(K_\Lambda).$$

Now, we have the following homeomorphisms of coset spaces : 
  $$ K_\Gamma/\mathbb{Z}^d \simeq  K_\Gamma/\big(K_\Gamma\cap\tilde{\Gamma}\big)\simeq K_\Gamma/\tilde{\Gamma} \subset \tilde{E'}/\tilde{\Gamma}$$
Since $\tilde{\Gamma}$ is a co-compact lattice in $\tilde{E'}$,  we have  $\tilde{E'}/\tilde{\Gamma}$  compact. Therefore, $K_\Gamma/\mathbb{Z}^d$ is a compact Lie subgroup of $E'$. Similarly,  $K_\Lambda/\mathbb{Z}^d$ is a compact Lie subgroup of $E'$. Since any compact subgroup of  a connected nilpotent group is a central torus ( Lemma 2.2.3, \cite{BFT16}), we obtain that $K_\Gamma/\mathbb{Z}^d$ and $K_\Lambda/\mathbb{Z}^d$ are both central torus. We consider the maximal torus 
$\mathbb{T}^d$ in $E'$. Since $\mbox{dim}\big( K_\Gamma)=\mbox{dim}(K_\Lambda)$, we have $K_\Gamma/\mathbb{Z}^d$ and $K_\Lambda/\mathbb{Z}^d$  two central torus inside the maximal torus $\mathbb{T}^d$  with same dimension.  Therefore, they are isomorphic as Lie groups, which implies that there is an isomorphism 
$$ K_\Gamma\simeq K_\Lambda.$$  
 Therefore, there is a Lie group isomorphism 
  $$  N(\Gamma)\longrightarrow N(\Lambda). $$
 Hence we prove our theorem. \hfill\(\square\)
 
 \section{Further Remarks}
 
 \subsection{Box spaces of nilpotent groups}
 
 The author proves in \cite{Das18} (Theorem 1.1) that if $\Gamma$ and $\Lambda$ are two finitely generated residually finite groups and a box space $\square_{\Gamma_n} \Gamma$ of $\Gamma$ and a box space $\square_{\Lambda_n} \Lambda$ of $\Lambda$ are Coarsely Equivalent, then $\Gamma$ and $\Lambda$ are UME. It was asked to the author by Y. de Cornulier when the converse is true. By Theorem \ref{imece}, we obtain that if $\Gamma$ and $\Lambda$ are finitely generated torsion free nilpotent groups which are IME (particularly, if they are UME), then any two box spaces  $\square_{\Gamma_n} \Gamma$ of $\Gamma$ and  $\square_{\Lambda_n} \Lambda$ of $\Lambda$ are Coarsely Equivalent. Therefore, we obtain the following theorem: 
 
 \begin{thm}
 Let $\Gamma$ and $\Lambda$ be two finitely generated torsion free nilpotent groups. Then the following are equivalent: 
  \begin{itemize}
 \item[1.] $\Gamma$ and $\Lambda$ are IME; 
 \item[2.] Any two box spaces  $\square_{\Gamma_n} \Gamma$ of $\Gamma$ and  $\square_{\Lambda_n} \Lambda$ of $\Lambda$ are Coarsely Equivalent;
 \item[3.] $\Gamma$ and $\Lambda$ are UME; 
 \end{itemize}
 \end{thm}
 
 \subsection{ QI embedding between two nilpotent groups}
 
It is well-known that  if there exists a continuous homo-morphism from a finitely generated torsion free nilpotent group $\Gamma$ to another such group $\Lambda$, then there is a Lie group homo-morphism from the Mal'cev completion of $N(\Gamma)$ of $\Gamma$ to Mal'cev completion $N(\Lambda)$ of $\Lambda$ (\cite{Mal51}. We can generalize this result under weaker assumption of  QI-embedding. 
 
 \begin{thm}
 If there exists an QI-embedding from a finitely generated torsion free nilpotent group $\Gamma$ to another such group $\Lambda$, then there is a Lie group homo-morphism from the Mal'cev completion of $N(\Gamma)$ of $\Gamma$ to Mal'cev completion $N(\Lambda)$ of $\Lambda$ 
 \end{thm}

 \subsection{Betti numbers of nilpotent groups}
 
 Using our main theorem \ref{MainTheorem} and Theorem 4.1.1 and Theorem 1.5 in \cite{Sha04}, we obtain the 
 following theorem
 
 \begin{thm}
  Let $\Gamma$ and $\Lambda$ be two finitely generated torsion free nilpotent groups. If $\Gamma$ and $\Lambda$ are IME, then the betti numbers
  $$b_n(\Gamma)=b_n(\Lambda)$$
  for all $n\in\mathbb{N}$ and the cohomological dimension with respect to co-effcients in a commutative divisible ring $R$ with unit
  $$\mbox{cd}_{R}(\Gamma)=\mbox{cd}_{R}(\Lambda).$$
 \end{thm}

 \vspace{5mm}

\textbf{Acknowledgement.} The author would like to thank  Romain Tessera for introducing  the article \cite{BFT16} and  Uri Bader for referring \cite{Gla03} and communicating \cite{Ba} with the author.

\end{document}